\documentclass[journal]{IEEEtran}
\ifCLASSINFOpdf
\else
\fi
\usepackage[cmex10]{amsmath}
\usepackage{amssymb}
\usepackage{amsfonts}
\usepackage{nicefrac}
\usepackage{mathrsfs}    
\usepackage{dsfont}
\usepackage{amsthm}
\usepackage{graphicx}

\newcommand{\normzero}[1]{\lVert #1 \rVert_{\ell_0}}

\providecommand{\abs}[1]{\lvert#1\rvert}

\newcommand{\bone}{\mathbf{1}}

\newcommand{\RP}{\bbR_+}
\newcommand{\bbR}{{\mathbb R}}

\newcommand{\ZP}{\mathbb{Z}_+}

\newcommand{\beps}{\varepsilon}

\newcommand{\crX}{\mathscr{X}}
\newcommand{\crY}{\mathscr{Y}}
\newcommand{\crC}{\mathscr{C}}

\newcommand{\set}[1]{\left\{ #1 \right\}}

\DeclareMathOperator*{\id}{id}

\DeclareMathOperator*{\supp}{supp}

\newtheorem{theorem}{Theorem}[section]
\newtheorem{lemma}{Lemma}[section]
\newtheorem{definition}{Definition}[section]
\newtheorem{example}{Example}[subsection]

\newtheorem{condition}{Condition}

\newcommand{\dist}{f}
\newcommand{\permgrp}{S_n}

\newcommand{\sparsity}{K}

\newcommand{\Supp}[1]{\supp\left(#1\right)}

\newcommand{\dime}{D_{\lambda}}
\newcommand{\rep}{M^{\lambda}}
\newcommand{\disthat}{\hat{\dist}(\lambda)}
\newcommand{\prop}[1]{\text{\sf P}(\sparsity,\lambda)}
\newcommand{\comp}{\crC}
\newcommand{\rand}{R(\sparsity, \comp)}
\newcommand{\rT}{R(\sparsity, T)}

\newcommand{\event}{\mathscr{E}}
\newcommand{\eveF}{\mathscr{F}}

\newcommand{\Prob}[1]{\Pr \left(#1\right)}

\newcommand{\R}{\mathbb{R}}

\begin{document}
\title{Inferring Rankings Using Constrained Sensing}
\author{Srikanth~Jagabathula and~Devavrat~Shah
  \thanks{This work was supported in parts by NSF CAREER CNS 0546590
    and NSF CCF 0728554.}
  \thanks{ Both SJ and DS are with the Department of
    Electrical Engineering and Computer Science, Massachusetts
    Institute of Technology, Cambridge, MA, 02140 USA e-mail:{\tt
      jskanth@alum.mit.edu, devavrat@mit.edu}}
}

\markboth{IEEE Transactions on Information Theory}%
{}
%



\maketitle

\begin{abstract}
  We consider the problem of recovering a function over the space of
  permutations (or, the symmetric group) over $n$ elements from given partial
  information; the partial information we consider is related to the group
  theoretic Fourier Transform of the function. This problem naturally arises in
  several settings such as ranked elections, multi-object tracking, ranking
  systems, and recommendation systems. Inspired by the work of Donoho and Stark
  in the context of discrete-time functions, we focus on non-negative functions
  with a sparse support (support size $\ll$ domain size). Our recovery method is
  based on finding the sparsest solution (through $\ell_0$ optimization) that is
  consistent with the available information. As the main result, we derive
  sufficient conditions for functions that can be recovered exactly from partial
  information through $\ell_0$ optimization. Under a natural random model for
  the generation of functions, we quantify the recoverability conditions by
  deriving bounds on the sparsity (support size) for which the function
  satisfies the sufficient conditions with a high probability as $n \to
  \infty$. $\ell_0$ optimization is computationally hard. Therefore, the popular
  compressive sensing literature considers solving the convex relaxation,
  $\ell_1$ optimization, to find the sparsest solution. However, we show that
  $\ell_1$ optimization fails to recover a function (even with constant
  sparsity) generated using the random model with a high probability as $n \to
  \infty$. In order to overcome this problem, we propose a novel iterative
  algorithm for the recovery of functions that satisfy the sufficient
  conditions. Finally, using an Information Theoretic framework, we study
  necessary conditions for exact recovery to be possible.
\end{abstract}

%
\IEEEpeerreviewmaketitle

\section{Introduction}
%
%
%
%
\IEEEPARstart{F}{unctions} over permutations serve as rich tools for modeling uncertainty in several important practical applications; they correspond to a general model class, where each model has a factorial number of parameters. However, in many practical applications, only partial information is available about the underlying functions; this is because either the problem setting naturally makes only partial information available, or memory constraints allow only partial information to be maintained as opposed to the entire function -- which requires storing a factorial number of parameters in general. In either case, the following important question arises: which ``types'' of functions can be recovered from access to only partial information? Intuitively, one expects a characterization that relates the ``complexity'' of the functions that can be recovered to the ``amount'' of partial information one has access to. One of the main goals of this paper is to formalize this statement. More specifically, this paper considers the problem of {\em exact} recovery of a function over the space of permutations given only partial information. When the function is a probability distribution, the partial information we consider can be thought of as lower-order marginals; more generally, the types of partial information we consider are related to the group theoretic Fourier Transform of the function, which provides a general way yo represent varying ``amounts'' of partial information. In this context, our goal is to (a) characterize a class of functions that can be recovered exactly from the given partial information, and (b) design a procedure for their recovery. We restrict ourselves to non-negative functions, which span many of the useful practical applications. Due to the generality of the setting we consider, a thorough understanding of this problem impacts a wide-ranging set of applications. Before we present the precise problem formulation and give an overview of our approach, we provide below a few motivating applications that can be modeled effectively using functions over permutations.

A popular application where functions over permutations naturally arise is the problem of {\em rank aggregation}. This problem arises in various contexts. The classical setting is that of {\em ranked election}, which has been studied in the area of {\em Social Choice Theory} for the past several decades. In the ranked election problem, the goal is to determine a ``socially preferred'' ranking of $n$ candidates contesting an election using the individual preference lists (permutations of candidates) of the voters. Since the ``socially preferred'' outcome should be independent of the identities of voters, the available information can be summarized as a function over permutations that maps each permutation $\sigma$ to the fraction of voters that have the preference list $\sigma$. While described in the context of elections, the ranked election setting is more general and also applies to aggregating through polls the population preferences on global issues, movies, movie stars, etc. Similarly, rank aggregation has also been studied in the context of aggregating webpage rankings~\cite{Dwork01}, where one has to aggregate rankings over a large number of webpages.  Bulk of the work done on the ranked election problem deals with the question of aggregation {\em given} access to the entire function over permutations that summarizes population preferences. In many practical settings, however, determining the function itself is non-trivial -- even for reasonable small values of $n$. Like in the setting of polling, one typically can gather only partial information about population preferences. Therefore, our ability to recover functions over permutations from available partial information impacts our ability to aggregate rankings. Interestingly, in the context of ranked election, Diaconis~\cite{D88} showed through spectral analysis that a partial set of Fourier coefficients of the function possesses ``rich'' information about the underlying function. This hints to the possibility that, in relevant applications, limited partial information can still capture a lot of structure of the underlying function. 

Another important problem, which has received a lot of attention recently, is the {\em Identity Management Problem} or the {\em Multi-object tracking problem}. This problem is motivated by applications in air traffic control and sensor networks, where the goal is to track the identities of $n$ objects from noisy measurements of identities and positions. Specifically, consider an area with sensors deployed that can identify the unique signature and the position associated with each object when it passes close to it. Let the objects be labeled $1, 2, \dotsc, n$ and let $x(t) = (x_1(t), x_2(t), \dotsc, x_n(t))$ denote the vector of positions of the $n$ objects at time $t$. Whenever a sensor registers the signature of an object the vector $x(t)$ is updated. A problem, however, arises when two objects, say $i, j$, pass close to a sensor simultaneously. Because the sensors are inexpensive, they tend to confuse the signatures of the two objects; thus, after the two objects pass, the sensor has information about the positions of the objects, but it only has beliefs about which position belongs to which object. This problem is typically modeled as a probability distribution over permutations, where, given a position vector $x(t)$, a permutation $\sigma$ of $1, 2, \dotsc, n$ describes the assignment of the positions to objects. Because the measurements are noisy, to each position vector $x(t)$, we assign, not a single permutation, but a distribution over permutations. Since we now have a distribution over permutations, the factorial blow-up makes it challenging to maintain it. Thus, it is often approximated using a partial set of Fourier coefficients. Recent work by~\cite{HGG08, KHJ07} deals with updating the distribution with new observations in the Fourier domain. In order to obtain the final beliefs one has to recover the distribution over permutations from a partial set of Fourier coefficients. 

Finally, consider the task of coming up with rankings for teams in a sports league, for example, the ``Formula-one'' car racing or American football, given the outcomes of various games. In this context, one approach is to model the final ranking of the teams using, not just one permutation, but a distribution over permutations. A similar approach has been taken in ranking players in online games (cf. Microsoft's TrueSkill solution \cite{trueskill}), where the authors, instead of maintaining scores, maintain a distribution over scores for each player. In this context, clearly, we can gather only partial information and the goal is to fit a model to this partial information. Similar questions arise in recommendation systems in cases where rankings, instead of ratings, are available or are preferred.

In summary, all the examples discussed above relate to inferring a function over permutations using partial information. To fix ideas, let $\permgrp$ denote the permutation group of order $n$ and $\dist \colon \permgrp \to \RP$ denote a non-negative function defined over the permutations. We assume we have access to partial information about $f(\cdot)$ that, as discussed subsequently, corresponds to a subset of coefficients of the group theoretic Fourier Transform of $f(\cdot)$. We note here that a partial set of Fourier coefficients not only provides a rigorous way to compress the high-dimensional function $f(\cdot)$ (as used in \cite{HGG08, KHJ07}), but also have natural interpretations, which makes it easy to gather in practice. Under this setup, our goal is to characterize the functions $\dist$ that can be recovered.  The problem of exact recovery of functions from a partial information has been widely studied in the context of discrete-time functions; however, the existing approaches don’t naturally extend to our setup. One of the classical approaches for recovery is to find the function with the minimum ``energy'' consistent with the given partial information. This approach was extended to functions over permutations in~\cite{KORT01}, where the authors obtain lower bounds on the energy contained in subsets of Fourier Transform coefficients to obtain better $\ell_2$ guarantees when using the function the minimum ``energy.'' This approach, however, does not naturally extend to the case of exact recovery. In another approach, which recently gained immense popularity, the function is assumed to have a sparse support and conditions are derived for the size of the support for which exact recovery is possible. This work was pioneered by Donoho; in~\cite{DS89}, Donoho and Stark use generalized uncertainty principles to recover a discrete-time function with sparse support from a limited set of Fourier coefficients. Inspired by this, we restrict our attention to functions with a sparse support.

Assuming that the function is sparse, our approach to performing exact recovery is to find the function with the sparsest support that is consistent with the given partial information, henceforth referred to as $\ell_0$ optimization. This approach is often justified by the philosophy of {\em Occam's razor}. We derive sufficient conditions in terms of sparsity (support size) for functions that can be recovered through $\ell_0$ optimization. Furthermore, finding a function with the sparsest support through $\ell_0$ minimization is in general computationally hard. This problem is typically overcome by considering the convex relaxation of the $\ell_0$ optimization problem. However, as we show in Theorem~\ref{thm:l1min}, such a convex relaxation does not yield exact recovery in our case. Thus, we propose a simple iterative algorithm called the `sparsest-fit' algorithm and prove that the algorithm performs exact recovery of functions that satisfy the sufficient conditions.

It is worth noting that our work has important connections to the work done in the recently popular area of {\em compressive sensing}. Broadly speaking, this work derives sufficient conditions under which the sparsest function that is consistent with the given information can be found by solving the corresponding $\ell_1$ relaxation problem. However, as discussed below in the section on relevant work, the sufficient conditions derived in this work do not apply to our setting. Therefore, our work may be viewed as presenting an alternate set of conditions under which the $\ell_0$ optimization problem can be solved efficiently. 

\subsection{Related Work}

Fitting sparse models to observed data has been a classical approach used in statistics for model recovery and is inspired by the philosophy of {\em Occam's Razor}. Motivated by this, sufficient conditions based on sparsity for learnability have been of great interest over years in the context of communication, signal processing and statistics, cf. \cite{S49, N02}. In recent years, this approach has become of particular interest due to exciting developments and wide ranging applications including: 
\begin{itemize} 
\item In signal processing (see \cite{CT05, CRT06, CR06, CRT06R, D06}) where the goal is to estimate a `signal' by means of minimal number of measurements. This is referred to as compressive sensing.  
\item In coding theory through the design of low-density parity check codes~\cite{G62, SS96, LMSS01} or in the design Reed Solomon codes \cite{RS60} where the aim is to design a coding scheme with maximal communication rate.  
\item In the context of streaming algorithms through the design of `sketches' (see \cite{T06, T04, BGIK08, CM06, GSTV07}) for the purpose of maintaining a minimal `memory state' for the streaming algorithm's operation.  
\end{itemize}

In all of the above work, the basic question (see \cite{M05}) pertains to the design of an $m \times n$ ``measurement'' matrix $A$ so that $x$ can be recovered efficiently from measurements $y = Ax$ (or its noisy version) using the ``fewest'' possible number measurements $m$. The setup of interest is when $x$ is sparse and when $m < n$ or $m \ll n$.  The type of interesting results (such as those cited above) pertain to characterization of the sparsity $K$ of $x$ that can be recovered for a given number of measurements $m$. The usual tension is between the ability to recover $x$ with large $k$ using a sensing matrix $A$ with minimal $m$.

The sparsest recovery approach of this paper is similar (in flavor) to the above stated work; in fact, as is shown subsequently, the partial information we consider can be written as a linear transform of the function $f(\cdot)$. However, the methods or approaches of the prior work do not apply. Specifically, the work considers finding the sparsest function consistent with the given partial information by solving the corresponding $\ell_1$ relaxation problem. The work derives a necessary and sufficient condition, called the {\em Restricted Nullspace Property}, on the structure of the matrix $A$ that guarantees that the solutions to the $\ell_0$ and $\ell_1$ relaxation problems are the same (see \cite{CRT06, BGIK08}). However, such sufficient conditions trivially fail in our setup (see ~\cite{JS08}). Therefore, our work provides an alternate set of conditions that guarantee efficient recovery of the sparsest function.

\subsection{Our Contributions}
Recovery of a function over permutations from only partial information is clearly a hard problem both from a theoretical and computational standpoint. We make several contributions in this paper to advance our understanding of the problem in both these respects.  As the main result, we obtain sufficient conditions -- in terms of sparsity -- for functions that can be recovered exactly from partial information. Specifically, our result establishes a relation between the ``complexity'' (as measured in sparsity) of the function that can be recovered and the ``amount'' of partial information available.

Our recovery scheme consists of finding the sparsest solution consistent with the given partial information through $\ell_0$ optimization. We derive sufficient conditions under which a function can be recovered through $\ell_0$ optimization. First, we state the sufficient conditions for recovery through $\ell_0$ optimization in terms of the structural properties of the functions. To understand the strength of the sufficient conditions, we propose a random generative model for functions with a given support size; we then obtain bounds on the size of the support for which a function generated according to the random generative model satisfies the sufficient conditions with a high probability. To our surprise, it is indeed possible to recover, with high probability, functions with seemingly large sparsity for given partial information (see precise statement of Theorems \ref{thm:case1}-\ref{thm:case4} for details).

Finding the sparsest solution through $\ell_0$ optimization is computationally hard. This problem is typically overcome by considering the $\ell_1$ convex relaxation of the $\ell_0$ optimization problem. However, as we show in Example~\ref{ex1}, $\ell_1$ relaxation does not always result in exact recovery, even when the the sparsity of the underlying function is only $4$. In fact, a necessary and sufficient condition for $\ell_1$ relaxation to yield the sparsest solution $x$ that satisfies the constraints $y = Ax$ is the so called Restricted Nullspace Condition (RNC) on the measurement matrix $A$; interestingly, the more popular Restricted Isoperimetric Property (RIP) on the measurement matrix $A$ is a sufficient condition. However, as shown below, the types of partial information we consider can be written as a linear transform of $f(\cdot)$. Therefore, Example~\ref{ex1} shows that in our setting, the measurement matrix does not satisfy RNC. It is natural to wonder if Example~\ref{ex1} is anomalous. We show that this is indeed not the case. Specifically, we show in Theorem~\ref{thm:l1min} that, with a high probability, $\ell_1$ relaxation fails to recover a function generated according to the random generative model. 

Since convex relaxations fail in recovery, we exploit the structural property of permutations to design a simple iterative algorithm called the `sparsest-fit' algorithm to perform recovery. We prove that the algorithm recovers a function from a partial set of its Fourier coefficients as long as the function satisfies the sufficient conditions. 

We also study the limitation of {\em any} recovery algorithm to recover a function exactly from a given form of partial information. Through an application of classical information theoretic Fano's inequality, we obtain a bound on the sparsity beyond which recovery is not {\em asymptotically reliable}; a recovery scheme is called asymptotically reliable if the probability of error asymptotically goes to $0$.

In summary, we obtain an intuitive characterization of the ``complexity'' (as measured in sparsity) of the functions that can be recovered from the given partial information. We show how $\ell_1$ relaxation fails in recovery in this setting. Hence, the sufficient conditions we derive correspond to an alternate set of conditions that guarantee efficient recovery of the sparsest function. 

\subsection{Organization}

Section \ref{sec:model} introduces the model, useful notations and the precise formulation of the problem. In Section \ref{sec:mainresults}, we provide the statements of our results. Section \ref{sec:algo} describes our iterative algorithm that can recover $\dist$ from $\disthat$ when certain conditions (see Condition \ref{cond1}) are satisfied.  Sections \ref{sec:thmzero} to \ref{sec:thmconverse} provide detailed proofs. Conclusions are presented Section \ref{sec:conclusion}.

\section{Problem Statement} \label{sec:model} 

In this section, we introduce the necessary notations, definitions and provide the formal problem statement.

\subsection{Notations} Let $n$ be the number of elements and $\permgrp$ be set of all possible $n!$ permutations or rankings of these of $n$ elements. Our interest is in learning non-negative valued functions $\dist$ defined on $S_n$, i.e. $\dist: \permgrp \to \RP$, where $\RP = \{x \in {\mathbb R} : x \geq 0 \}$. The support of $\dist$ is defined as
\begin{equation*}
\Supp{\dist} = \set{\sigma \in \permgrp \colon \dist(\sigma) \neq 0}.  
\end{equation*}
The cardinality of support, $|\Supp{\dist}|$ will be called the {\em sparsity} of $\dist$ and will be denoted by $K$. We will also call it the $\ell_0$ norm of $f$, denoted by $\lvert f \rvert_0$.

In this paper, we wish to learn $\dist$ from a partial set of Fourier coefficients. To define the Fourier transform of a function over the permutation group, we need some notations.  To this end, consider a partition of $n$, i.e. an ordered tuple $\lambda = (\lambda_1, \lambda_2, \ldots, \lambda_r)$, such that $\lambda_1 \geq \lambda_2 \geq \ldots \geq \lambda_r \geq 1$, and $n = \lambda_1 + \lambda_2 + \ldots + \lambda_r$. For example, $\lambda = (n-1,1)$ is a partition of $n$. Now consider a partition of the $n$ elements, $\{1,\dots,n\}$, as per the $\lambda$ partition, i.e. divide $n$ elements into $r$ bins with $i$th bin having $\lambda_i$ elements. It is easy to see that $n$ elements can be divided as per the $\lambda$ partition in $D_\lambda$ distinct ways, with
$$ D_\lambda = \frac{n!}{ \prod_{i=1}^r \lambda_i !}. $$
Let the distinct partitions be denoted by $t_i, 1\leq i\leq D_\lambda$\footnote{To keep notation
simple, we use $t_i$ instead of $t_i^\lambda$ that takes explicit 
dependence on $\lambda$ into account.}.  
For example, for $\lambda = (n-1,1)$ there are $D_\lambda = n!/(n-1)! = n$ 
distinct ways given by  
$$t_i \equiv \{1,\dots,i-1,i+1,\dots,n\}\{i\}, ~1\leq i\leq n.$$
Given a permutation $\sigma \in \permgrp$, its action on $t_i$ is defined through its action on the $n$ elements of $t_i$, resulting in a $\lambda$ partition with the $n$ elements permuted.  In the above example with $\lambda = (n-1,1)$, $\sigma$ acts on $t_i$ to give $t_{\sigma(i)}$, i.e.  
\begin{align*}
  \sigma: t_i \to t_{\sigma(i)}, \text{ where } t_i \equiv \{1,\dots,i-1,i+1,\dots,n\}\{i\} \text{ and } \\
 t_{\sigma(i)} \equiv \{1,\dots,\sigma(i)-1,\sigma(i)+1,\dots,n\}\{\sigma(i)\}.
\end{align*}
Now, for a given partition $\lambda$ and a permutation $\sigma \in \permgrp$, define a $0/1$ valued $D_\lambda \times D_\lambda$ matrix $\rep(\sigma)$ as \begin{equation*}
  \rep_{ij}(\sigma) = \begin{cases} 1,& \text{if }\sigma(t_j) = t_i \\
0,&\text{otherwise.}
\end{cases} 
~~ ~~\text{for all } 1 \leq i,j \leq \dime
\end{equation*}
This matrix $\rep(\sigma)$ corresponds to a degree $D_\lambda$ representation of the permutation group.

\subsection{Partial Information as a Fourier Coefficient} The partial information we consider in this paper is the Fourier transform coefficient of $\dist$ at the representation $\rep$, for each $\lambda$. The motivation for considering Fourier coefficients at representations $\rep$ is two fold: first, they provide a rigorous way to compress the high-dimensional function $f(\cdot)$ (as used in \cite{HGG08, KHJ07}), and second, as we shall see, Fourier coefficients at representations $\rep$ have natural interpretations, which makes it easy to gather in practice. In addition, each representation $\rep$ contains a subset of the lower-order irreducible representations; thus, for each $\lambda$, $\rep$ conveniently captures the information contained in a subset of the lower-order Fourier coefficients up to $\lambda$. We now define the Fourier coefficient of $\dist$ at the representation $\rep$, which we call $\lambda$-partial information. 

\begin{definition}[$\lambda$-Partial Information]
Given a function $\dist \colon \permgrp \to \RP$ and partition $\lambda$. The Fourier Transform coefficient at representation $\rep$, which we call  the $\lambda$-partial information, is denoted by $\disthat$ and is defined as 
$$ \disthat = \sum_{\sigma \in \permgrp} f(\sigma) \rep(\sigma). $$
\end{definition}

Recall the example of $\lambda = (n-1,1)$ with $f$ as a probability distribution on $\permgrp$. Then, $\disthat$ is an $n\times n$ matrix with the $(i,j)$th entry being the probability of element $j$ mapped to element $i$ under $f$. That is, $\disthat$ corresponds to the {\em first order} marginal of $f$ in this case.

\subsection{Problem Formulation}\label{prob:form}

We wish to recover a function $\dist$ based on its partial information $\disthat$ based on partition $\lambda$. As noted earlier, the classical approach based on Occam's razor suggests recovering the function as a solution of the following $\ell_0$ optimization problem: \begin{eqnarray}\label{l0opt}
  {\sf minimize} & & \|g\|_0 \qquad {\sf over} \qquad g : \permgrp \to \RP \nonumber \\
  {\sf subject ~to} & & \hat{g}(\lambda) = \disthat.  \end{eqnarray}

We note that the question of recovering $\dist$ from $\disthat$ is very similar to the question studied in the context of compressed sensing, i.e. recover $x$ from $y = Ax$. To see this, with an abuse of notation imagine $\disthat$ as the $D_\lambda^2$ dimensional vector and $\dist$ as $n!$ dimensional vector. Then, $\disthat = A \dist$ where each column of $A$ corresponds to $\rep(\sigma)$ for certain permutation $\sigma$. The key difference from the compressed sensing literature is that $A$ is given in our setup rather than being a design choice.

\vspace{.1in} {\em Question One.} ~As the first question of interest, we wish to identify precise conditions under which $\ell_0$ optimization problem \eqref{l0opt} recovers the original function $f$ as its unique solution.

Unlike the popular literature (cf. compressed sensing), such conditions can not be based on sparsity only. This is well explained by the following (counter-)example. In addition, the example also shows that linear independence of the support of $f$ does not guarantee uniqueness of the solution to the $\ell_0$ optimization problem.

\begin{example}\label{ex1}
  For any $n \geq 4$, consider the four permutations $\sigma_1 = (1,2)$, $\sigma_2 = (3,4)$, $\sigma_3 = (1,2)(3,4)$ and $\sigma_4 = \id$, where $\id$ is the identity permutation.  In addition, consider the partition $\lambda = (n-1, 1)$.  Then, it is easy to see that
$$\rep(\sigma_1) + \rep(\sigma_2) = \rep(\sigma_3) + \rep(\sigma_4).$$ 

We now consider three cases where a bound on sparsity is not sufficient to guarantee the existence of a unique solution to \eqref{l0opt}.

\begin{enumerate}
\item This example shows that a sparsity bound (even $4$) on $f$ is not sufficient to guarantee that $f$ will indeed be the sparsest solution. Specifically, suppose that $\dist(\sigma_i) = p_i,$ where $p_i \in \RP$ for $1 \leq i \leq 4$, and $\dist(\sigma) = 0$ for all other $\sigma \in \permgrp$. Without loss of generality, let $p_1 \leq p_2$. 
Then, \begin{align*}
  &\disthat \\
= &p_1 \rep(\sigma_1) + p_2 \rep(\sigma_2) + p_3
  \rep(\sigma_3) + p_4 \rep(\sigma_4)\\
  =& (p_2 - p_1) \rep(\sigma_2) + (p_3 + p_1) \rep(\sigma_3) \\
&+ (p_4 + p_1) \rep(\sigma_4). 
\end{align*}
Thus, function $g$ with $g(\sigma_2) = p_2 - p_1$, $g(\sigma_3) = p_3+p_1$, $g(\sigma_4) = p_4 + p_1$ and $g(\sigma) = 0$ for all other $\sigma \in \permgrp$ is such that $\hat{g}(\lambda) = \disthat$ but $\|g\|_0 = 3 < 4 = \|\dist\|_0$. That is, $\dist$ can not be recovered as the solution of $\ell_0$ optimization problem \eqref{l0opt} even when support of $\dist$ is only $4$.
\item This example shows that although $f$ might be a sparsest solution, it may not be unique. In particular, suppose that $\dist(\sigma_1) = \dist(\sigma_2) = p$ and $\dist(\sigma) = 0$ for all other $\sigma \in \permgrp$. Then, $\disthat = p \rep(\sigma_1) + p \rep(\sigma_2) = p \rep(\sigma_3) + p \rep(\sigma_4)$. Thus, \eqref{l0opt} does not have a unique solution. 
\item Finally, this example shows that even though the support of $f$ corresponds to a linearly independent set of columns, the sparsest solution may not be unique. Now suppose that $\dist(\sigma_i) = p_i,$ where $p_i \in \RP$ for $1 \leq i \leq 3$, and $\dist(\sigma) = 0$ for all other $\sigma \in \permgrp$. Without loss of generality, let $p_1 \leq p_2$. Then,
\begin{align*}
  &\disthat \\
= &p_1 \rep(\sigma_1) + p_2 \rep(\sigma_2) + p_3
  \rep(\sigma_3) \\
 = &(p_2 - p_1) \rep(\sigma_2) + (p_3 + p_1) \rep(\sigma_3) + 
  p_1 \rep(\sigma_4). 
\end{align*}
Here, note that $\set{\rep(\sigma_1), \rep(\sigma_2), \rep(\sigma_3)}$ is linearly independent, yet the solution to \eqref{l0opt} is not unique. 
\end{enumerate}
\end{example}

\vspace{.1in} {\em Question Two.} ~The resolution of the first question will provide a way to recover $\dist$ by means of solving the $\ell_0$ optimization problem in \eqref{l0opt}. However, in general, it is computationally a hard problem. Therefore, we wish to obtain a simple and possibly iterative algorithm to recover $\dist$ (and hence for solving \eqref{l0opt}).

\vspace{.1in} {\em Question Three.} Once we identify the conditions for exact recovery of $\dist$, the next natural question to ask is ``how restrictive are the conditions we imposed on $\dist$ for exact recovery?'' In other words, as mentioned above, we know that the sufficient conditions don't translate to a simple sparsity bound on functions, however, can we find a sparsity bound such that ``most,'' if not all, functions that satisfy the sparsity bound can be recovered? We make the notion of ``most'' functions precise by proposing a natural random generative model for functions with a given sparsity. Then, for given a partition $\lambda$, we want to obtain $K(\lambda)$ so that if $K < K(\lambda)$ then recovery of $\dist$ generated according to the generative model from $\disthat$ is possible with high probability.

\vspace{.1in} This question is essentially an inquiry into whether the situation demonstrated by Example~\ref{ex1} is contrived or not. In other words, it is an inquiry into whether such examples happen with vanishingly low probability for a randomly chosen function. To this end, we describe a natural random function generation model.  
\begin{definition}[Random Model]\label{def:rmodel} Given $\sparsity \in \ZP$
  and an interval $\comp = [a,b],~0 < a < b$, a random function 
  $\dist$ with sparsity $\sparsity$ and values in $\comp$ is generated
  as follows: choose $\sparsity$ permutations from $\permgrp$ independently 
  and uniformly at random
  \footnote{Throughout, we will assume that 
  the random selection is done  {\em with} replacement. 
  }, say $\sigma_1,\dots, \sigma_\sparsity$; 
  select $\sparsity$ values from $\comp$ uniformly at random,
  say $p_1,\dots,p_\sparsity$; then function $f$ is defined as 
  $$ f(\sigma) = \begin{cases} p_i  & \text{if} \quad \sigma = \sigma_i, ~1\leq i\leq \sparsity \\
                               0  & \text{otherwise.}
                               \end{cases} $$
We will denote this model as $\rand$. 
\end{definition}

\vspace{.1in} {\em Question Four.} Can we characterize a limitation on the ability of {\em any} algorithm to recover $\dist$ from $\disthat$ ?

\section{Main Results} \label{sec:mainresults}

As the main result of this paper, we provide answers to the four questions stated in Section \ref{prob:form}.  We start with recalling some notations.  Let $\lambda = (\lambda_1,\dots,\lambda_r)$ be the given partition of $n$. We wish to recover function $\dist : S_n \to \R_+$ from available information $\disthat$.  Let the sparsity of $f$ be $\sparsity$,
$$\Supp{\dist} = \{\sigma_1,\dots, \sigma_\sparsity\}, \quad\mbox{and}\quad 
f(\sigma_k) = p_k, 1\leq k \leq \sparsity.$$

\vspace{.1in} {\em Answers One \& Two.} To answer the first two questions, we need to find sufficiency conditions for recovering $f$ through $\ell_0$ optimization \eqref{l0opt} and a simple algorithm to recover the function. For that, we first try to gain a qualitative understanding of the conditions that $f$ must satisfy. Note that a necessary condition for $\ell_0$ optimization to recover $f$ is that \eqref{l0opt} must have a {\em unique} solution; otherwise, without any additional information, we wouldn't know which of the multiple solutions is the true solution. Clearly, since $\disthat = \sum_{\sigma \in \permgrp} \dist(\sigma) \rep(\sigma)$, \eqref{l0opt} will have a unique solution only if $\set{\rep(\sigma)}_{\sigma \in \Supp{\dist}}$ is linearly independent. However, this linear independence condition is, in general, not sufficient to guarantee a unique solution; in particular, even if $\set{\rep(\sigma)}_{\sigma \in \Supp{\dist}}$ is linearly independent, there could exist $\set{\rep(\sigma')}_{\sigma' \in \mathcal{H}}$ such that $\disthat = \sum_{\sigma' \in \mathcal{H}} \rep(\sigma')$ and $\lvert \mathcal{H} \rvert \leq \sparsity$, where $\sparsity := \lvert \Supp{\dist} \rvert$; Example \ref{ex1} illustrates such a scenario. Thus, a sufficient condition for $\dist$ to be the unique sparsest solution of \eqref{l0opt} is that not only is $\set{\rep(\sigma)}_{\sigma \in \Supp{\dist}}$ linearly independent, but $\set{\rep(\sigma), \rep(\sigma')}_{\sigma \in \Supp{\dist}, \sigma' \in \mathcal{H}}$ is linearly independent for all $\mathcal{H} \subset \permgrp$ such that $\lvert \mathcal{H} \rvert \leq \sparsity$; in other words, not only we want $\rep(\sigma)$ for $\sigma \in \Supp{\dist}$ to be linearly independent, but we want them to be linearly independent even after the addition of at most $\sparsity$ permutations to the support of $\dist$. Note that this condition is similar to the Restricted Isometry Property (RIP) introduced in \cite{CT05}, which roughly translates to the property that $\ell_0$ optimization recovers $x$ of sparsity $K$ from $y = Ax$ provided every subset of $2K$ columns of $A$ is linearly independent. Motivated by this, we impose the following conditions on $\dist$.

\begin{condition}[Sufficiency Conditions]\label{cond1}{\em 
Let $f$ satisfy the following:
\begin{itemize}

\item[$\circ$] {\em Unique Witness}: for any $\sigma \in \Supp{\dist}$,
there exists $1\leq i_\sigma, j_\sigma \leq D_\lambda$ such that 
$M^\lambda_{i_\sigma j_\sigma}(\sigma) = 1$, but $ 
M^\lambda_{i_\sigma j_\sigma}(\sigma') = 0$, for all $\sigma' (\neq \sigma) \in \Supp{\dist}.$

\item[$\circ$] {\em Linear Independence}: for any collection of integers 
$c_1,\dots, c_\sparsity$ taking values in $\{-\sparsity, \dots, \sparsity\}$, 
$\sum_{k=1}^\sparsity c_k p_k \neq 0$, unless $c_1 =\ldots = c_\sparsity = 0$. 

\end{itemize}
}
\end{condition}

The above discussion motivates the ``unique witness'' condition; indeed, $\rep(\sigma)$ for $\sigma$ satisfying the ``unique witness'' condition are linearly independent because every permutation has a unique witness and no non-zero linear combination of $\rep(\sigma)$ can yield zero. On the other hand, as shown in the proof of Theorem \ref{thm:zero}, the {\em linear independence} condition is required for the uniqueness of the sparsest solution. 


Now we state a formal result that establishes Condition \ref{cond1} as sufficient for recovery of $f$ as the unique solution of $\ell_0$ optimization problem. Further, it allows for a simple, iterative recovery algorithm. Thus, Theorem \ref{thm:zero} provides answers to questions {\em One} and {\em Two} of Section \ref{prob:form}.

\begin{theorem}\label{thm:zero} 
Under Condition \ref{cond1}, the function $f$ is the unique solution of the $\ell_0$ optimization problem \eqref{l0opt}. Further, a simple, iterative algorithm called the sparsest-fit algorithm, described in Section \ref{sec:algo}, recovers $f$.  \end{theorem}

\vspace{.1in} 
{\em Linear Programs Don't Work.} Theorem \ref{thm:zero} states that under Condition \ref{cond1}, the $\ell_0$ optimization recovers $f$ and the  sparsest-fit algorithm is a simple iterative algorithm to recover it. In the context of compressive sensing literature (cf. \cite{CRT06,CRT06R,D06,BGIK08}), it has been shown that convex relaxation of $\ell_0$ optimization, such as the Linear Programing relaxation, have the same solution in similar scenarios.  Therefore, it is natural to wonder whether such a relaxation would work in our case. To this end, consider the following Linear Programing relaxation of \eqref{l0opt} stated as the following $\ell_1$ minimization problem: \begin{eqnarray}\label{l1opt}
  {\sf minimize} & & \|g\|_1 \qquad {\sf over} \qquad g : \permgrp \to \RP \nonumber \\
  {\sf subject ~to} & & \hat{g}(\lambda) = \disthat.  \end{eqnarray} 
Example~\ref{ex1} provides a scenario where $\ell_1$ relaxation fails in recovery. In fact, we can prove a stronger result. The following result establishes that -- with a high probability -- a function generated randomly as per Definition~\ref{def:rmodel} cannot be recovered by solving the linear program~\eqref{l1opt} because there exists a function $g$ such that $\hat{g}(\lambda) = \disthat$ and $\| g \|_1 = \| f \|_1$.

\begin{theorem} \label{thm:l1min} 
Consider a function $\dist$ randomly generated as per Definition \ref{def:rmodel} with sparsity $\sparsity \geq 2$. Then, as longs as $\lambda$ is not the partition  $(1, 1, \dotsc, 1)$ ($n$ times), with probability $1 - o(1)$, there exists a function $g$ distinct from $f$ such that $\hat{g}(\lambda) = \disthat$ and $\| g \|_1  = \| f \|_1$.
\end{theorem}

\vspace{.1in} {\em Answer Three.} Next, we turn to the third question.  Specifically, we study the conditions for high probability recoverability of a random function $f$ in terms of its sparsity. That is, we wish to identify the high probability recoverability threshold $K(\lambda)$. In what follows, we spell out the result starting with few specific cases so as to better explain the dependency of $K(\lambda)$ on $D_\lambda$.

\vspace{.1in} {\em Case 1}: $\lambda = (n-1,1)$. Here $D_\lambda = n$ and $\disthat$ provides the {\em first order} marginal information. As stated next, for this case the achievable recoverability threshold $K(\lambda)$ scales\footnote{Throughout this paper, by $\log $ we mean the natural logarithm, i.e. $\log_e$, unless otherwise stated.} as $n \log n$.

\begin{theorem}\label{thm:case1}
  A randomly generated $f$ as per Definition \ref{def:rmodel} can be recovered by the sparsest-fit algorithm with probability $1-o(1)$ as long as $\sparsity \leq (1-\beps) n \log n$ for any fixed $\beps > 0$.
\end{theorem}

\vspace{.1in} {\em Case 2}: $\lambda = (n-m,m)$ with $1 < m = O(1)$. Here $D_\lambda = \Theta(n^m)$ and $\disthat$ provides the {\em $m$th order} marginal information.  As stated next, for this case we find that $K(\lambda)$ scales at least as $n^m \log n$.

\begin{theorem}\label{thm:case2} 
A randomly generated $f$ as per Definition \ref{def:rmodel} can be recovered from $\disthat$ by the sparsest-fit algorithm for $\lambda = (n-m, m), m = O(1)$, with probability $1-o(1)$ as long as $\sparsity \leq \frac{(1-\beps)}{m!} n^m \log n$ for any fixed $\beps > 0$.  
\end{theorem}

In general, for any $\lambda$ with $\lambda_1 = n-m$ and $m=O(1)$, arguments of Theorem \ref{thm:case2} can be adapted to show that $K(\lambda)$ scales as $n^m \log n$.  Theorems \ref{thm:case1} and \ref{thm:case2} suggest that the recoverability threshold scales $D_\lambda \log D_\lambda$ for $\lambda = (\lambda_1,\dots,\lambda_r)$ with $\lambda_1 = n -m$ for $m = O(1)$. Next, we consider the case of more general $\lambda$.

\vspace{.1in}
{\em Case 3}: $\lambda = (\lambda_1,\dots,\lambda_r)$ with $\lambda_1 = n - O\left(n^{\frac{2}{9}-\delta}\right)$
for any $\delta > 0$. As stated next, for this case, the recoverability threshold 
$K(\lambda)$ scales at least as $D_\lambda \log\log D_\lambda$.

\begin{theorem}\label{thm:case3} 
A randomly generated $f$ as per Definition \ref{def:rmodel} can be recovered from $\disthat$ by the sparsest-fit algorithm for $\lambda = (\lambda_1,\dots,\lambda_r)$ with $ \lambda_1 = n - n^{\frac{2}{9}-\delta}$ for any $\delta > 0$, with probability $1-o(1)$ as long as $\sparsity \leq (1-\beps) D_\lambda \log\log D_\lambda$ for any fixed $\beps > 0$.  
\end{theorem}

\vspace{.1in} {\em Case 4}: Any $\lambda = (\lambda_1,\dots,\lambda_r)$. The results stated thus far suggest that the threshold is essentially $D_\lambda$, ignoring the logarithm term. For general $\lambda$, we establish a bound on $K(\lambda)$ as stated in Theorem \ref{thm:case4} below. Before stating the result, we introduce some notation. For given $\lambda$, define $\alpha = (\alpha_1,\dots, \alpha_r)$ with $\alpha_i = \lambda_i/n, ~1\leq i\leq r$. Let
$$ H(\alpha) = -\sum_{i=1}^r \alpha_i \log \alpha_i, \quad \text{and}\quad 
H'(\alpha) = - \sum_{i=2}^r \alpha_i \log \alpha_i.$$

\begin{theorem}\label{thm:case4} 
Given $\lambda = (\lambda_1,\dots,\lambda_r)$, a randomly generated $f$ as per Definition \ref{def:rmodel} can be recovered from $\disthat$ by the sparsest-fit algorithm with probability $1-o(1)$ as long as \begin{eqnarray}
 \sparsity & \leq & C\,D_\lambda^{\gamma(\alpha)},  \label{eq:thmcase4}
\end{eqnarray}
where
\begin{align*}
  \gamma(\alpha)  &=  \frac{M}{M+1} \left[ 1 - C' \frac{H(\alpha) - H'(\alpha)}{H(\alpha)} \right], 
\\ &\text{with } M = \left\lfloor\frac{1}{1-\alpha_1}\right\rfloor
\end{align*}
and $0 < C, C' < \infty$ are constants.
\end{theorem}

At a first glance, the above result seems very different from the crisp formulas of Theorems \ref{thm:case1}-\ref{thm:case3}.  Therefore, let us consider a few special cases. First, observe that as $\alpha_1 \uparrow 1$, $M/(M+1) \to 1$. Further, as stated in Lemma \ref{lem:analysis}, $H'(\alpha)/H(\alpha) \to 1$. Thus, we find that the bound on sparsity essentially scales as $D_\lambda$. Note that the cases 1, 2 and 3 fall squarely under this scenario since $\alpha_1 = \lambda_1/n = 1 - o(1)$. Thus, this general result contains the results of Theorems \ref{thm:case1}-\ref{thm:case3} (ignoring the logarithm terms).

Next, consider the other extreme of $\alpha_1 \downarrow 0$. Then, $M \to 1$ and again by Lemma \ref{lem:analysis}, $H'(\alpha)/H(\alpha) \to 1$.  Therefore, the bound on sparsity scales as $\sqrt{D_\lambda}$. This ought to be the case because for $\lambda = (1,\dots,1)$ we have $\alpha_1 = 1/n \to 1$, $D_\lambda = n!$, and unique witness property holds only up to $o(\sqrt{D_\lambda}) = o(\sqrt{n!})$ due to the standard Birthday paradox.

In summary, Theorem \ref{thm:case4} appears reasonably tight for the general form of partial information $\lambda$. We now state the Lemma \ref{lem:analysis} used above (proof in Appendix \ref{ap1}).
\begin{lemma}\label{lem:analysis}
  Consider any $\alpha = (\alpha_1,\dots,\alpha_r)$ with $1\geq \alpha_1\geq \dots \geq \alpha_r \geq 0$ and $\sum_{i=1}^r \alpha_r =1$. Then,
$$\lim_{\alpha_1 \uparrow 1} \frac{H'(\alpha)}{H(\alpha)} = 1, $$
$$\lim_{\alpha_1 \downarrow 0} \frac{H'(\alpha)}{H(\alpha)} = 1. $$
\end{lemma}

\vspace{.1in} 
{\em Answer Four.} Finally, we wish to understand the fundamental limitation
on the ability to recover $\dist$ from $\disthat$ by any algorithm. To 
obtain a meaningful bound (cf. Example \ref{ex1}), we shall examine this 
question under an appropriate information theoretic setup.  

To this end, as in random model $\rand$, consider a function 
$\dist$ generated with given $\sparsity$ and $\lambda$. For
technical reasons (or limitations), we will assume that the
values $p_i$s are chosen from a discrete set. Specifically, 
let each $p_i$ be chosen from integers $\{1,\dots, T\}$ instead 
of compact set $\comp$. We will denote this random model by 
$\rT$. 

Consider any algorithm that attempts to recover $\dist$ 
from $\disthat$ under $\rT$. Let $h$ be the estimation 
of the algorithm. Define probability of error of the
algorithm as 
$$ p_{\text{err}} = \Pr\left(h \neq \dist\right).$$
We state the following result. 
\begin{theorem}\label{thm:converse}
With respect to random model $\rT$, the probability of
error is uniformly bounded away from $0$ for all $n$ large
enough and any $\lambda$, if 
\begin{eqnarray}
\sparsity & \geq & \frac{3 D_{\lambda}^2}{ n \log n} \left[ \log \left(\frac{D_{\lambda}^2}{ n \log n} \vee T\right) \right], \end{eqnarray}
where for any two numbers $x$ and $y$, $x \vee y$ denotes $\max\set{x, y}$. 
\end{theorem}

\section{Sparsest-fit algorithm}\label{sec:algo}
As mentioned above, finding the sparsest distribution that is
consistent with the given partial information is in general a
computationally hard problem. In this section, we propose an efficient
algorithm to fit the sparsest distribution to the given partial
information $\hat{f}(\lambda)$, for any partition $\lambda$ of $n$. 
The algorithm we propose determines the sparsest distribution {\em
  exactly} as long as the underlying distribution belongs to the
general family of distributions that satisfy the `unique witness' and
`linear independence' conditions; we call this the `sparsest-fit'
algorithm. In this case, it follows from Theorem~\ref{thm:zero} that
the `sparsest-fit' algorithm indeed recovers the underlying
distribution $f(\cdot)$ exactly from partial information
$\hat{f}(\lambda)$. 
When the conditions are not satisfied, the algorithm produces a
certificate to that effect and aborts.  

Using the degree $D_\lambda$ representation of the permutations, the algorithm
processes the elements of the partial information matrix $\hat{f}(\lambda)$
sequentially and incrementally builds the permutations in the support. We
describe the sparsest-fit algorithm as a general procedure to recover a set of
non-negative values given sums of these values over a collection of subsets,
which for brevity we call subset sums. In this sense, it can be thought of as a
linear equation solver customized for a special class of systems of linear
equations.

Next we describe the algorithm in detail and prove the relavant
theorems.

\subsection{Sparsest-fit algorithm} \label{sec:sparsest-fit}
We now describe the sparsest-fit algorithm that was also referred to in Theorems
\ref{thm:zero}, \ref{thm:case1}-\ref{thm:case4} to recover function $\dist$ from
$\disthat$ under Condition~\ref{cond1}.

\vspace{.1in} {\em Setup.} The formal description of the algorithm is given in
Fig.~\ref{fig:algo}. The algorithm is described there as a generic procedure
to recover a set of non-negative values given a collection of their subset sums.
As explained in Fig.~\ref{fig:algo}, the inputs to the algorithm are $L$
positive numbers $q_1,\dots, q_L$ sorted in ascending order $q_1 \leq q_2 \leq
\dotsb \leq q_L$. As stated in assumptions C1-C3 in Fig.~\ref{fig:algo}, the
algorithm assumes that the $L$ numbers are different subset sums of $\sparsity$
distinct positive numbers $p_1,\dots, p_\sparsity$ i.e., $q_\ell = \sum_{T_\ell}
p_k$ for some $T_\ell \subset \set{1, 2, \dotsc,
  K}$, and the values and subsets satisfy the conditions: for each $1 \leq k
\leq K$, $p_k = q_\ell$ for some $1 \leq \ell \leq L$ and $\sum_{T} p_k \neq
\sum_{T'} p_k$ for $T \neq T'$. Given this setup, the sparsest-fit algorithm
recovers the values $p_k$ and subset membership sets $A_k := \set{\ell \colon k
  \in T_\ell}$ for $1 \leq k \leq \sparsity$ using $q_\ell$, but without any
knowledge of $\sparsity$ or subsets $T_\ell, 1\leq \ell \leq L$.

Before we describe the algorithm, note that in order to use the sparsest-fit
algorithm to recover $f(\cdot)$ we give the non-zero elements of the partial
information matrix $\hat{f}(\lambda)$ as inputs $q_\ell$. In this case, $L$
equals the number of non-zero entries of $\hat{f}(\lambda)$, $p_k =
f(\sigma_k)$, and the sets $A_k$ correspond to $\rep(\sigma_k)$. Here,
assumption C1 of the algorithm is trivially satisfied. As we argue in
Section~\ref{sec:thmzero}, assumptions C2, C3 are implied by the `unique
witness' and `linear independence' conditions.

\vspace{.1in} {\em Description.} The formal description is given below in the
Fig.~\ref{fig:algo}. The algorithm processes elements $q_1, q_2, \dotsc, q_L$
sequentially and builds membership sets incrementally. It maintains the number
of non-empty membership sets at the end of each iteration $\ell$ as $k(\ell)$.
Partial membership sets are maintained as sets $A_k$, which at the end of
iteration $\ell$ equals $\set{1 \leq k \leq k(\ell) \colon k \in T_{\ell'}
  \text{ for some } \ell' \leq \ell}$. The values found are maintained as $p_1,
p_2, \dotsc, p_{k(\ell)}$. The value of $k(0)$ is initialized to zero and the
sets $A_k$ are initialized to be empty.

In each iteration $\ell$, the algorithm checks if the value $q_\ell$ can be
written as a subset sum of values $p_1, p_2, \dotsc, p_{k(\ell - 1)}$ for some
subset $T$. If $q_\ell$ can be expressed as $\sum_{k \in T} p_k$ for some $T
\subset \set{1, 2, \dotsc, k(\ell -1)}$, then the algorithm adds $\ell$ to sets
$A_k$ for $k \in T$ and updates $k(\ell)$ as $k(\ell) = k(\ell - 1)$ before
ending the iteration. In case there exists no such subset $T$, the algorithm
updates $k(\ell)$ as $k(\ell -1) + 1$, makes the set $A_{k(\ell)}$ non-empty by
adding $\ell$ to it, and sets $p_{k(\ell)}$ to $q_\ell$. At the end the
algorithm outputs $(p_k, A_k)$ for $1 \leq k \leq k(L)$.

\begin{figure}[!h]
\vspace{0.2in}
\vspace{3pt}
\noindent\framebox[\linewidth][t]{
  \begin{minipage}{0.97\linewidth} 
{\noindent \bf Input:} Positive values
    $\set{q_1, q_2,
  \ldots, q_L}$ sorted in ascending order i.e., $q_1\leq q_2 \leq \dotsc \leq q_L$.
\vspace{5pt}

{\noindent \bf Assumptions:} 
$\exists$ positive values $\set{p_1, p_2, \ldots, p_\sparsity}$ such that:
\begin{enumerate}
\item[C1.] For each $1 \leq \ell \leq L$, $q_\ell = \sum_{k \in T_\ell} p_k$, for
  some $T_\ell \subseteq \set{1, 2, \ldots, \sparsity}$
\item[C2.] For each $1 \leq k \leq \sparsity$, there exists a $q_\ell$ such that $q_\ell
  = p_\sparsity$.
\item[C3.] $\sum_{k \in T} p_k \neq \sum_{k' \in T'} p_{k'}$, for all $T, T'
  \subseteq \set{1, 2, \ldots, J}$ and $T \cap T' = \emptyset$.
\end{enumerate}

\vspace{12pt}

{\noindent \bf Output:} $\set{p_1, p_2, \ldots, p_\sparsity}$, $\forall ~ 1 \leq k \leq \sparsity$ set $A_k$ s.t.
$$A_k = \{ \ell : q_\ell = \sum_{j \in T} p_j \text{ and index } k
\text{ belongs to set } T\}.$$ 

\vspace{12pt}

\noindent{\bf Algorithm:} \\


{\em initialization}: $p_0 = 0, k(0) = 0$, $A_k = \emptyset$ for all possible $k$. \\
{\bf for} ~~ $\ell = 1$ ~to~ $L$ \\
\hspace*{0.5cm} {\bf if} ~ $q_\ell = \sum_{k \in T} p_k$ ~for~ some~ $T \subseteq \{0, 1, \ldots, k(\ell-1) \}$ \\
\hspace*{1cm} $k(\ell) = k(\ell-1)$ \\
\hspace*{1cm} $A_k = A_k \cup \set{\ell}$ ~~$\forall$~~ $k \in T$ \\
\hspace*{0.5cm} {\bf else} \\
\hspace*{1cm} $k(\ell) = k(\ell-1) + 1$ \\
\hspace*{1cm} $p_{k(\ell)} = q_\ell$ \\
\hspace*{1cm} $A_{k(\ell)} = A_{k(\ell)} \cup \set{\ell}$ \\
\hspace*{0.5cm} {\bf end if} \\
{\bf end for} \\
{\em Output} $\sparsity = k(L)$ and $(p_k, A_k), 1\leq k \leq \sparsity$.
\vspace{3pt}
\end{minipage}}
\caption{Sparsest-fit algorithm}
\label{fig:algo}
\end{figure}

We now argue that under assumptions C1-C3 stated in Fig.~\ref{fig:algo}, the
algorithm finds $(p_k, A_k)$ for $1 \leq k \leq K$ accurately. Note that by
Assumption C2, there exists at least one $q_\ell$ such that it is equal to
$p_k$, for each $1 \leq k \leq \sparsity$. Assumption C3 guarantees that the
condition in the {\bf if} statement is not satisfied whenever $q_\ell =
p_{k(\ell)}$. Therefore, the algorithm correctly assigns values to each of the
$p_k$s. Note that the condition in the {\bf if} statement being true implies
that $q_\ell$ is a subset sum of some subset $T \subset \set{p_1, p_2, \ldots,
  p_{k(\ell-1)}}$. Assumption C3 ensures that if such a combination exists then
it is unique. Thus, when the condition is satisfied, index $\ell$ belongs only
to the sets $A_k$ such that $k \in T$. When the condition in the {\bf if}
statement is false, then from Assumptions C2 and C3 it follows that $\ell$ is
contained only in $A_{k(\ell)}$. From this discussion we conclude that the
sparsest-fit algorithm correctly assigns all the indices to each of the $A_k$s.
Thus, the algorithm recovers $p_k, A_k$ for $1\leq k\leq \sparsity$ under
Assumptions C1, C2 and C3. We summarize it in the following Lemma.

\begin{lemma}\label{lem:algo} 
  The sparsest-fit algorithm recovers $p_k, A_k$ for
  $1\leq k\leq \sparsity$ under Assumptions C1, C2 and C3. 
\end{lemma}

{\em Complexity of the algorithm.} Initially, we sort at most $D_{\lambda}^2$
elements. This has a complexity of $O(D_{\lambda}^2 \log D_{\lambda})$. Further,
note that the {\bf for} loop in the algorithm iterates for at most
$D_{\lambda}^2$ times. In each iteration, we are solving a subset-sum problem.
Since there are at most $\sparsity$ elements, the worst-case complexity of
subset-sum in each iteration is $O(2^{\sparsity})$. Thus, the worst-case
complexity of the algorithm is $O( D_{\lambda}^2 \log D_{\lambda} +
D_{\lambda}^2 2^{\sparsity})$. However, using the standard balls and bins
argument, we can prove that for $K = O(D_\lambda \log D_\lambda)$, with a high
probability, there are at most $O(\log D_\lambda)$ elements in each subset-sum
problem. Thus, the complexity would then be $O\left(\exp(\log^2
  D_\lambda)\right)$ with a high probability.

\section{Proof of Theorem \ref{thm:zero}}\label{sec:thmzero}

The proof of Theorem \ref{thm:zero} requires us to establish two
claims : under Condition \ref{cond1}, (i) the sparsest-fit algorithm finds $\dist$ and (ii) the $\ell_0$ optimization \eqref{l0opt} has $\dist$ as it's unique solution. We establish these two claims in that order.

\vspace{.1in} {\em The sparsest-fit algorithm works.} As noted in Section \ref{sec:algo}, the sparsest-fit algorithm can be used to recover $\dist$ from $\disthat$. As per Lemma \ref{lem:algo}, the correctness of the sparsest-fit algorithm follows under Assumptions C1, C2 and C3. The Assumption C1 is trivially satisfied in the context of recovering $\dist$ from $\disthat$ as discussed in Section \ref{sec:algo}. Next, we show that Condition \ref{cond1} implies C2 and C3.  Note that the {\em unique witness} of Condition \ref{cond1} implies C2 while C3 is a direct implication of {\em linear independence} of Condition \ref{cond1}. Therefore, we have established that the sparsest-fit algorithm recovers $\dist$ from $\disthat$ under Condition \ref{cond1}.

\vspace{.1in}
{\em Unique Solution of $\ell_0$ Optimization.}  To arrive at a 
contradiction, assume that there exists a function 
$g  \colon \permgrp \to \RP$ such that 
$\hat{g}(\lambda) = \disthat$ and $ L \stackrel{\triangle}{=} 
\normzero{g} \leq \normzero{f} = \sparsity$. 
Let 
$$\Supp{\dist} = \{\sigma_k \in \permgrp : 1\leq k \leq \sparsity\}, 
~\dist(\sigma_k) = p_k, 1\leq k\leq \sparsity,$$ 
$$ \Supp{g} = \{\rho_\ell \in \permgrp : 1\leq \ell \leq L \}, ~ g(\rho_\ell) = q_\ell, 1\leq \ell\leq L.$$
By hypothesis of Theorem \ref{thm:zero}, $\dist$ satisfies Condition \ref{cond1}. 
Therefore, entries of matrix $\disthat$ contains the values 
$p_1, p_2, \ldots, p_{\sparsity}$. Also, by our assumption 
$\disthat = \hat{g}(\lambda)$. Now, by definition, each entry of the
matrix $\hat{g}(\lambda)$ is a summation of a subset of $L$
numbers, $q_\ell, 1\leq \ell\leq L$. Therefore, it follows that 
for each $k, 1\leq k\leq \sparsity$, we have
$$ p_k = \sum_{j \in T_k} q_j, \qquad \text{for some} \quad T_k \subseteq \set{1, 2, \ldots,
    L}.$$ 
Equivalently, 
\begin{eqnarray} 
p & = & A q,\label{eq:zero1}
\end{eqnarray}
 where $p = [p_k]_{1\leq k\leq \sparsity}$, $q = [q_\ell]_{1\leq \ell\leq L}$ 
$A \in \{0,1\}^{\sparsity \times L}$. 

Now consider the matrix $\disthat$. As noted before, each of its entries is a summation of a subset of numbers $p_k, 1\leq k\leq \sparsity$. Further, each $p_k, 1\leq k\leq \sparsity$ contributes to exactly $D_\lambda$ distinct entries of $\disthat$. Therefore, it follows that the summation of all entries of $\disthat$ is $D_\lambda (p_1+\dots+p_\sparsity)$. That is,
$$ \sum_{ij} \disthat_{ij} ~=~ D_\lambda \left(\sum_{k=1}^\sparsity p_k\right).$$
Similarly, 
$$ \sum_{ij} \hat{g}(\lambda)_{ij} ~=~ D_\lambda \left(\sum_{\ell=1}^L q_\ell\right).$$
But $\disthat = \hat{g}(\lambda)$. Therefore, 
\begin{eqnarray}
p \cdot \bone & = & q \cdot \bone, \label{eq:zero2}
\end{eqnarray}
where $\bone$ is vector of all $1$s of appropriate dimension (we have
abused the notation $\bone$ here): in LHS, it is of dimension $\sparsity$, 
in RHS it is of dimension $L$. Also, from \eqref{eq:zero1} we have
\begin{eqnarray}
p \cdot \bone & = & Aq \cdot \bone \nonumber \\
              & = & \sum_{\ell=1}^L c_\ell q_\ell, \label{eq:zero3}
\end{eqnarray}
for some $c_j \in \ZP$. From \eqref{eq:zero2} and \eqref{eq:zero3},
it follows that 
\begin{eqnarray}
\sum_j q_j & = & \sum_j c_j q_j.\label{eq:zero4}
\end{eqnarray}
  Now, there are two options: (1) either all the $c_\ell$s are $> 0$, 
  or (2) some of them are equal to zero. 
  In the case (1), when $c_\ell > 0$ for all $1 \leq \ell \leq L$, 
  it follows that $c_\ell = 1$ for each 
  $1\leq \ell \leq L$; or else, RHS of \eqref{eq:zero4} will be strictly
  larger than LHS since $q_\ell > 0$ for all $1\leq \ell\leq L$ by
  definition. Therefore, the matrix $A$ in \eqref{eq:zero1} must contain
  exactly one non-zero entry, i.e. $1$, in each column. Since
  $p_k > 0$ for all $1\leq k\leq \sparsity$, it follows that there must
  be at least $\sparsity$ non-zero entries in $A$. Finally, since
  $L \leq \sparsity$, it follows that we must have $L = \sparsity$. 
  In summary, it must be that $A$ is a $\sparsity \times \sparsity$
  matrix with each row and column having exactly one $1$, and rest of the
  entries $0$. That is, $A$ is a permutation matrix. That is, 
  $p_k, 1\leq k\leq \sparsity$ is permutation of $q_1,\dots, q_L$ with
  $L = \sparsity$. By relabeling the $q_\ell$s, if required, without 
  loss of generality, we assume that $p_k = q_k$, for $1 \leq k \leq \sparsity$. 
  Since $\hat{g}(\lambda) = \disthat$ and $p_k = q_k$ for $1 \leq k \leq \sparsity$, 
  it follows that $g$ also satisfies Condition 
  \ref{cond1}. Therefore, the sparsest-fit algorithm accurately recovers $g$ 
  from $\hat{g}(\lambda)$. Since the input to the algorithm is only 
  $\hat{g}(\lambda)$ and $\hat{g}(\lambda) = \disthat$, it follows that 
  $g = f$ and we have reached contradiction to our assumption that $f$
  is not the unique solution of optimization problem \eqref{l0opt}. 
  
  Now consider the remaining case (2) and suppose that $c_\ell = 0$
  for some $\ell$. Then, it follows that some of the columns in the 
  $A$ matrix are zeros. Removing those columns of $A$ we can 
  write 
  $$ p = \tilde{A} \tilde{q}, $$ 
  where $\tilde{A}$ is formed from $A$ by removing the zero columns 
  and $\tilde{q}$ is formed from $q$ by removing $q_\ell$s such that 
  $c_\ell = 0$. Let $\tilde{L}$ be the size of $\tilde{q}$. Since at 
  least one column was removed, $\tilde{L} < L \leq \sparsity$. The 
  condition $\tilde{L} < \sparsity$ implies that the vector $p$ lies 
  in a lower dimensional space. Further, $\tilde{A}$ is a $0, 1$
  valued matrix. Therefore, it follows that $p$ violates 
  the linear independence property of Condition \ref{cond1} resulting in a contradiction. This completes the proof of Theorem 
  \ref{thm:zero}. 

\section{Proof of Theorem \ref{thm:l1min}}
We prove this theorem by showing that when two permutations, say $\sigma_1, \sigma_2$, are chosen uniformly at random, with a high probability, the sum of their representation matrices $\rep(\sigma_1) + \rep(\sigma_2)$ can be decomposed in at least two ways. For that, note that a permutation can be represented using cycle notation, e.g. for $n=4$, the permutation $1\mapsto 2, 2\mapsto 1, 3\mapsto 4, 4\mapsto 3$ can be represented as a composition of two cycles $(1 2) (3 4)$.  We call two cycles {\em distinct} if they have no elements in common, e.g. the cycles $(1 2)$ and $(3 4)$ are distinct. Given two permutations $\sigma_1$ and $\sigma_2$, let $\sigma_{1,2} = \sigma_1 \sigma_2$ be their composition.

Now consider two permutations $\sigma_1$ and $\sigma_2$ such 
that they have distinct cycles. For example, $\sigma_1 = (1,2)$ 
and $\sigma_2 = (3,4)$ are permutations with distinct
cycles. Then $\sigma_{1,2} = \sigma_1 \sigma_2 = (1 2) (3 4)$. 
We first prove the theorem for $\lambda = (n-1, 1)$ and then extend it to a general $\lambda$;
thus, we fix the partition $\lambda = (n-1, 1)$. Then, we have:
\begin{equation}
  \label{eq:4}
  \rep(\sigma_1) + \rep(\sigma_2) = \rep(\sigma_{1,2}) + \rep(\id)
\end{equation}
where $\sigma_1$ and $\sigma_2$ have distinct cycles and $\id$ is the
identity permutation. Now, assuming that $p_1 \leq p_2$, consider the
following:
\begin{align*}
  &p_1 \rep(\sigma_1) + p_2 \rep(\sigma_2) \\
   =~ &
  p_1 \rep(\sigma_{1,2}) + p_1 \rep(\id) + (p_2 - p_1) \rep(\sigma_2).
\end{align*}
Thus, given $\disthat = p_1 \rep(\sigma_1) + p_2 \rep(\sigma_2)$, it
can be decomposed in two distinct ways with both having the same
$\ell_1$ norm. Of course, the same analysis can be carried out when
$\dist$ has a sparsity $K$. Thus, we conclude that whenever $\dist$
has two permutations with distinct cycles in its support, the $\ell_1$
minimization solution is not unique.  Therefore, to establish claim 
of Theorem \ref{thm:l1min}, it is sufficient to prove that when 
we choose two permutations uniformly at random, they have distinct 
cycles with a high probability. 

To this end, let $\event$ denote the event that two permutations 
chosen uniformly at random have distinct cycles. Since permutations 
are chosen uniformly at random, $\Prob{\event}$ can be computed 
by fixing one of the permutations to be $\id$. Then, $\Prob{\event}$ 
is the probability that a permutation chosen at random has 
more than one cycle. 

Let us evaluate $\Prob{\event^c}$. For that, consider a permutation
having exactly one cycle with the cycle containing $l$ elements. The
number of such permutations will be $\binom{n}{l} (l-1)!$. This is
because we can choose the $l$ elements that form the cycle in
$\binom{n}{l}$ ways and the $l$ numbers can be arranged in the cycle
in $(l-1)!$ ways. Therefore,
\begin{equation}
  \label{eq:nips0818}
  \Pr(\event^c) = \frac{1}{n!} \sum_{l= 1}^n \binom{n}{l} (l-1)! = \sum_{r=1}^n \frac{1}{l (n-l)!}
\end{equation}
Now, without loss of generality let's assume that $n$ is even. Then,
\begin{equation}
  \label{eq:nips0819}
  \sum_{l=1}^{n/2} \frac{1}{l (n-l)!} \leq \sum_{l=1}^{n/2} \frac{1}{\left(\frac{n}{2}\right)!} = \frac{1}{\left(\frac{n}{2} - 1 \right)! }
\end{equation}
The other half of the sum becomes
\begin{equation}
  \label{eq:nips0820}
  \sum_{l=n/2}^n \frac{1}{l(n-l)!} \leq \sum_{k=0}^{n/2} \frac{1}{ \frac{n}{2}  k!}  \leq \frac{2}{n} \sum_{k=0}^{\infty} \frac{1}{k!} \leq \frac{O(1)}{n}
\end{equation}
Putting everything together, we have
\begin{align*}
  \Pr(\event) \geq  1 - \Pr(\event^c) &\geq 1 - 
\left( \frac{1}{\left( \frac{n}{2} -1 \right)! } + \frac{O(1)}{n}
\right) \\
&\to 1 \text{ as } n \to \infty.
\end{align*}
Thus, Theorem \ref{thm:l1min} is true for $\lambda = (n-1, 1)$.

In order to extend the proof to a general $\lambda$, we observe that the standard cycle notation for a permutation we discussed above can be extended to $\lambda$ partitions for a general $\lambda$. Specifically, for any given $\lambda$, observe that a permutation can be imagined as a perfect matching in a $D_{\lambda} \times D_{\lambda}$ bipartite graph, which we call the $\lambda$-bipartite graph and denote it by $G^{\lambda} = (V_1^{\lambda} \times V_2^{\lambda}, E^{\lambda})$; here $V_1^{\lambda}$ and $V_2^{\lambda}$ respectively denote the left and right vertex sets with $\lvert V_1^{\lambda} \rvert = \lvert V_2^{\lambda} \rvert = D_{\lambda}$ with a node for every $\lambda$ partition of $n$. Let $t_1, t_2, \dotsc, t_{D_{\lambda}}$ denote the $D_{\lambda}$ $\lambda$-partitions of $n$; then, the nodes in $V_1^{\lambda}$ and $V_2^{\lambda}$ can be labeled by $t_1, t_2, \dotsc, t_{D_{\lambda}}$. Since every perfect matching in a bipartite graph can be decomposed into its corresponding distinct cycles (the cycles can be obtained by superposing the bipartite graph corresponding to identity permutation with the $\lambda$-bipartite graph of the permutation), every permutation can be written as a combination of distinct cycles in its $\lambda$-bipartite graph. The special case of this for $\lambda = (n-1,1)$ is the standard cycle notation we discussed above; for brevity, we call the $\lambda$-bipartite graph for $\lambda = (n-1,1)$ the standard bipartite graph.

In order to prove the theorem for a general $\lambda$, using an argument similar to above, it can be shown that it is sufficient to prove that a randomly chosen permutation contains at least two distinct cycles in its $\lambda$-bipartite graph with a high probability. For that, it is sufficient to prove that a permutation with at least two distinct cycles in its standard bipartite graph has at least two distinct cycles in its $\lambda$-bipartite graph for any general $\lambda$. The theorem then follows from the result we established above that a randomly chosen permutation has at least two distinct cycles in its standard bipartite graph with a high probability.

To that end, consider a permutation, $\sigma$, with at least two distinct cycles in the standard bipartite graph. Let $A := (a_1, a_2, \dotsc, a_{\ell_1})$ and $B := (b_1, b_2, \dotsc, b_{\ell_2})$ denote the first two cycles in the standard bipartite graph; clearly, $\ell_1 \ell_2 \geq 2$ and at least one of $\ell_1, \ell_2$ is $\leq n/2$. Without loss of generality we assume that $\ell_2 \leq n/2$. Let $\lambda = (\lambda_1, \lambda_2, \dotsc, \lambda_r)$. Since $\lambda_1 \geq \lambda_2 \geq \ldots \geq \lambda_r$, we have $\lambda_r \leq n/2$. First, we consider the case when $\lambda_r < n/2$. Now consider the $\lambda$-partition, $t_1$, of $n$ constructed as follows: $a_1$ placed in the $r$th partition, $a_2$ in the first partition, all the elements of the second cycle $b_1, b_2, \dotsc, b_{\ell_2}$ arbitrarily in the first $r-1$ partitions and the rest placed arbitrarily. Note that such a construction is possible by the assumption on $\lambda_r$. Let $t_1'$ denote $\sigma(t_1)$; then, $t_1' \neq t_1$ because $t_1$ does not contain $a_2$ in the $r$th partition while $t_1'$ contains $\sigma(a_1) = a_2$ in the $r$th partition. Thus, the partition $t_1$ belongs to a cycle that has a length of at least $2$ partitions. Thus, we have found one cycle, which we denote by $C_1$. Now consider a second partition $t_2$ constructed as follows: $b_1$ placed in the $r$th partition, $b_2$ in the first and the rest placed arbitrarily. Again, note that $\sigma(t_2) \neq t_2$. Thus, $t_2$ belongs to a cycle of length at least $2$, which we denote by $C_2$. Now we have found two cycles $C_1, C_2$, and we are left with proving that they are distinct. In order to establish the cycles are distinct, note that none of the partitions in cycle $C_1$ can be $t_2$. This is true because, by construction, $t_2$ contains $b_1$ in the $r$th partition while none of the partitions in $C_1$ can contain any elements from the cycle $B$ in the $r$th partition. This finishes the proof for all $\lambda$ such that $\lambda_r < n/2$.

We now consider the case when $\lambda_r = n/2$. Since $\lambda_1 \geq
\lambda_r$, it follows that $r = 2$ and $\lambda = (n/2, n/2)$. For
$\ell_2 < n/2$, it is still feasible to construct $t_1$ and $t_2$, and
the theorem follows from the arguments above. Now we consider the case
when $\ell_1 = \ell_2 = n/2$; let $\ell := \ell_1 = \ell_2$. Note that
now it is infeasible to construct $t_1$ as described above. Therefore,
we consider $t_1 = \set{a_1, b_2, \dotsc, b_{\ell}}\set{b_1, a_2,
  \dotsc, a_{\ell}}$ and $t_2 = \set{b_1, a_2, \dotsc,
  a_{\ell}}\set{a_1, b_2, \dotsc, b_{\ell}}$. Clearly, $t_1 \neq t_2$,
$\sigma(t_1) \neq t_1$ and $\sigma(t_2) \neq t_2$. Thus, $t_1$ and
$t_2$ belong to two cycles, $C_1$ and $C_2$, each with length at least
$2$. It is easy to see that these cycles are also distinct because
every $\lambda-$partition in the cycle $C_1$ will have only one
element from cycle $A$ in the first partition and, hence, $C_1$ cannot
contain the $\lambda-$partition $t_2$. This completes the proof of the
theorem. 

\section{Proof of Theorem \ref{thm:case1} : $\lambda = (n-1,1)$} 
Our interest is in recovering a random function $\dist$ from 
partial information $\disthat$. To this end, let 
\begin{align*}
\sparsity = \|\dist\|_0, ~~ &\Supp{\dist} = \{\sigma_k \in \permgrp : 1\leq k \leq \sparsity\}, 
\\~~ \text{and} ~~&\dist(\sigma_k) = p_k, 1\leq k\leq \sparsity.  
\end{align*}
Here $\sigma_k$ and $p_k$ are randomly chosen as per the random
model $\rand$ described in Section \ref{sec:model}. For $\lambda = (n-1,1)$,
$D_\lambda = n$; then $\disthat$ is an $n\times n$ matrix with its
$(i,j)$th entry being 
$$ \disthat_{ij} = \sum_{k: \sigma_k(j) = i} p_k, \qquad \text{for}\qquad 1\leq i, j \leq n. $$

To establish Theorem \ref{thm:case1}, we prove that 
as long as $\sparsity \leq C_1 n \log n$ with $C_1 = 1 -\beps$, 
$\dist$ can be recovered by the sparsest-fit algorithm with probability $1-o(1)$ for any
fixed $\beps >0$.  
Specifically, we show that for $\sparsity \leq C_1 n\log n$, 
Condition \ref{cond1} is satisfied with probability $1-o(1)$, which in turn implies that  
the sparsest-fit algorithm recovers $\dist$ as per Theorem \ref{thm:zero}. 
Note that the ``linear independence'' property of Condition
\ref{cond1} is satisfied with probability $1$ under
$\rand$ as $p_k$ are chosen from a distribution 
with continuous support. Therefore, we are left with
establishing ``unique witness'' property. 

To this end, let $4\delta = \beps$ so  that $C_1 \leq 1-4\delta$.
 Let  $\event_k$ be
the event that $\sigma_k$ satisfies the unique witness
property, $1\leq k\leq \sparsity$. Under $\rand$, 
since $\sparsity$ permutations are chosen from $\permgrp$
independently and uniformly at random, it follows that 
$\Pr(\event_k)$ is the same for all $k$. Therefore, by union 
bound, it is sufficient to establish that 
$\sparsity \Pr(\event^c_1) = o(1)$. Since we are 
interested in $\sparsity = O(n\log n)$, it is sufficient 
to establish $\Pr(\event^c_1) = O(1/n^2)$.  Finally, once again 
due the symmetry, it is sufficient to evaluate $\Pr(\event_1)$
assuming $\sigma_1 = \id$, i.e. $\sigma_1(i) = i$ for all 
$1\leq i\leq n$. Define 
$$\eveF_j = \{\sigma_k(j)\neq j, \quad \text{for}\quad 2\leq k\leq \sparsity\}, ~~\text{for}~~1\leq j\leq n.$$
It then follows that 
$$ \Pr(\event_1) = \Pr\left(\cup_{j=1}^n \eveF_j\right).$$
Therefore, for any $L \leq n$, we have 
\begin{eqnarray} 
\Pr(\event_1^c) & = & \Pr\left(\cap_{j=1}^n \eveF_j^c\right) \nonumber \\
               & \leq & \Pr\left(\cap_{j=1}^L \eveF_j^c\right) \nonumber \\
              & = & \Pr\left(\eveF_1^c\right) \left[\prod_{j=2}^L \Pr\left(\eveF_{j}^c \,\Big|\, \cap_{\ell=1}^{j-1} \eveF_\ell^c\right) \right]. \label{eq:c11}
\end{eqnarray}
Next we show that for the selection of $L = n^{1-\delta}$, the RHS
of \eqref{eq:c11} is bounded above by $\exp(-n^\delta) = O(1/n^2)$. 
That will complete the proof of achievability. 

For that, we start by bounding $\Pr(\eveF_1^c)$:
\begin{eqnarray}
\Pr\left(\eveF_1^c\right) & = & 1 - \Pr\left(\eveF_1\right) \nonumber \\
& = & 1 - \left(1 - \frac{1}{n}\right)^{\sparsity-1}. \label{eq:c12}
\end{eqnarray}
The last equality follows because all permutations are chosen
uniformly at random. For $j \geq 2$, we now evaluate 
$\Pr\left(\eveF_{j}^c \,\Big|\, \cap_{\ell=1}^{j-1} \eveF_\ell^c\right)$. 
Given $\cap_{\ell=1}^{j-1} \eveF_\ell^c$, for any $k, 2\leq k\leq \sparsity$, 
$\sigma_k(j)$ will take a value from $n - j +1$ values, possibly including $j$, uniformly at random. Thus, we obtain the following
bound:
\begin{eqnarray}
\Pr\left(\eveF_{j}^c \,\Big|\, \cap_{\ell=1}^{j-1} \eveF_\ell^c\right)
& \leq & 1 - \left(1-\frac{1}{n-j+1}\right)^{\sparsity-1}. \label{eq:c13}
\end{eqnarray}
From \eqref{eq:c11}-\eqref{eq:c13}, we obtain that 
\begin{eqnarray} 
\Pr(\event_1^c) & \leq & \prod_{j=1}^{L} \left(1 - \left(1-\frac{1}{n-j+1}\right)^{\sparsity-1}\right) \nonumber \\
                & \leq  & \left[1-\left(1-\frac{1}{n-L}\right)^{\sparsity}\right]^L \nonumber \\
                & \leq & \left[1-\left(1-\frac{1}{n-L}\right)^{C_1 n \log n}\right]^L, \label{eq:c14}
\end{eqnarray}                
where we have used $\sparsity \leq C_1 n \log n$ in the last inequality. 
Since $L = n^{1-\delta}$, $n - L = n (1-o(1))$. Using the standard fact
$1-x = e^{-x} (1+O(x^2))$ for small $x \in [0,1)$, we have
\begin{eqnarray}
\left(1-\frac{1}{n-L}\right) & = & \exp\left(-\frac{1}{n-L}\right) \left(1 + O\left(\frac{1}{n^2}\right)\right). \label{eq:c15}
\end{eqnarray}
Finally, observe that 
$$ \left(1 + O\left(\frac{1}{n^2}\right)\right)^{C_1 n \log n} = \Theta(1).$$
Therefore, from \eqref{eq:c14} and \eqref{eq:c15}, it follows that 
\begin{eqnarray} 
\Pr(\event_1^c) & \leq & \left[1- \Theta\left(\exp\left(-\frac{C_1 \log n}{1-n^{-\delta}}\right)\right) \right]^L \nonumber \\
& \leq & \left[1 - \Theta\left(\exp\left(-(C_1+\delta) \log n\right)\right)\right]^L \nonumber \\
& = & \left[1 - \Theta\left(\frac{1}{n^{C_1+\delta}}\right)\right]^L \nonumber \\
& \leq & \exp\left(-\Theta\left(\frac{L}{n^{C_1 + \delta}}\right)\right) \nonumber \\
& = & \exp\left(-\Omega(n^{2\delta})\right), \label{eq:c16}
\end{eqnarray}
where we have used the fact that $1-x\leq e^{-x}$ for $x \in [0,1]$ and 
$L = n^{1-\delta}, ~C_1 \leq 1-4\delta$. From \eqref{eq:c16}, it follows that
$\Pr(\event_1) = O(1/n^2)$. This completes the proof of achievability. 

\section{Proof of Theorem \ref{thm:case2} : $\lambda = (n-m,m)$} 
Our interest is in recovering the random function $\dist$ from 
partial information $\disthat$. As in proof of Theorem \ref{thm:case1},
we use the notation  
\begin{align*}
\sparsity = \|\dist\|_0, ~~ &\Supp{\dist} = \{\sigma_k \in \permgrp :
1\leq k \leq \sparsity\}, \\
\text{and} ~~&\dist(\sigma_k) = p_k, 1\leq k\leq \sparsity.  
\end{align*}
Here $\sigma_k$ and $p_k$ are randomly chosen as per the random
model $\rand$ described in Section \ref{sec:model}. For $\lambda = (n-m,m)$,
$D_\lambda = \frac{n!}{(n-m)! m!} ~\sim~n^m$ and $\disthat$ is an 
$D_\lambda \times D_\lambda$ matrix. 

To establish Theorem \ref{thm:case2}, we shall prove that 
as long as $\sparsity \leq C_1 n^m \log n$ with $0 < C_1 < \frac{1}{m!}$
a constant, $\dist$ can be recovered by the sparsest-fit algorithm with probability $1-o(1)$.
We shall do so by verifying that the Condition \ref{cond1} holds
with probability $1-o(1)$, so that the sparsest-fit algorithm will recover $\dist$ 
as per Theorem \ref{thm:zero}. As noted earlier, the 
``linear independence'' of Condition \ref{cond1} is satisfied 
with probability $1$ under $\rand$. Therefore, we are left with 
establishing the ``unique witness'' property. 

To this end, for the purpose of bounding, without loss of
generality, let us assume that $\sparsity = \frac{(1-2\delta)}{m!}n^m \log n$
for some $\delta > 0$. Set $L = n^{1-\delta}$. 
Following arguments similar to those in the proof of Theorem \ref{thm:case1}, 
it will be sufficient to establish that $\Pr(\event_1^c = O(1/n^{2m})$; 
where $\event_1$ is the event that permutation $\sigma_1 = \id$
satisfies the unique witness property. 

To this end, recall that $\disthat$ is a $D_\lambda \times D_\lambda$
matrix. Each row (and column) of this matrix corresponds to a distinct
$\lambda$ partition of $n$ : $t_i, 1\leq i\leq D_\lambda$.
Without loss of generality, let us order
the $D_\lambda$ $\lambda$ partitions of $n$ so that the $i$th
partition, $t_i$, is defined as follows: $ t_1 = \{1,\dots,n-m\}\{n-m+1,\dots, n\}, $
and for $2\leq i\leq L$, 
\begin{multline*}
 t_i = \{1,\dots,n-im,n-(i-1)m+1,\dots,n\}\\\{n-im+1,\dotsc, n-(i-1)m\}.
\end{multline*}
Note that since $\sigma_1 = \id$, we have 
$\sigma_1(t_i) = t_i$ for all $1\leq i\leq D_\lambda$.
Define 
$$\eveF_j = \{\sigma_k(t_j)\neq t_j, \quad \text{for}\quad 2\leq k\leq \sparsity\}, ~~\text{for}~~1\leq j\leq D_\lambda.$$
Then it follows that 
$$ \Pr(\event_1) = \Pr\left(\cup_{j=1}^{D_\lambda} \eveF_j\right).$$
Therefore,  
\begin{eqnarray} 
\Pr(\event_1^c) & = & \Pr\left(\cap_{j=1}^{D_\lambda} \eveF_j^c\right) \nonumber \\
               & \leq & \Pr\left(\cap_{j=1}^L \eveF_j^c\right) \nonumber \\
              & = & \Pr\left(\eveF_1^c\right) \left[\prod_{j=2}^L \Pr\left(\eveF_{j}^c \,\Big|\, \cap_{\ell=1}^{j-1} \eveF_\ell^c\right) \right]. \label{eq:c21}
\end{eqnarray}
First, we bound $\Pr(\eveF_1^c)$. Each permutation $\sigma_k, k \neq 1$, 
maps $t_1 = \{1,\dots,n-m\}\{n-m+1,\dots,n\}$ to 
$\{\sigma_k(1),\dots\sigma_k(n-m)\}\{\sigma_k(n-m+1),\dots,\sigma_k(n)\}$. 
Therefore, $\sigma_k(t_1) = t_1$ iff $\sigma_k$ maps set of elements 
$\{n-m+1,\dots,n\}$ to the same set of elements. Therefore, 
\begin{eqnarray}
 \Pr\left(\sigma_k(t_1) = t_1\right) & = & \frac{1}{{n \choose m}} \nonumber \\
          & =  & \frac{m!}{\prod_{\ell=0}^{m-1} (n-\ell)}. \nonumber \\
          & \leq & \frac{m!}{(n-Lm)^m}.\label{eq:c22}
\end{eqnarray}          
Therefore, it follows that 
\begin{eqnarray}
\Pr\left(\eveF_1^c\right) & = & 1 - \Pr\left(\eveF_1\right) \nonumber \\
& = & 1 - \Pr\left(\sigma_k(t_1) \neq t_1, ~2\leq k\leq \sparsity\right) \nonumber \\
& = & 1 - \prod_{k=2}^{\sparsity} \left(1 - \Pr\left(\sigma_k(t_1) = t_1\right)\right) \nonumber \\
& \leq & 1 - \left(1 - \frac{m!}{(n-Lm)^m}\right)^{\sparsity}. \label{eq:c23}
\end{eqnarray}
Next we evaluate $\Pr\left(\eveF_{j}^c \,\Big|\, \cap_{\ell=1}^{j-1} \eveF_\ell^c\right)$
for $2\leq j\leq L$. Given $\cap_{\ell=1}^{j-1} \eveF_\ell^c$, we have (at least partial)
information about the action of $\sigma_k, 2\leq k\leq \sparsity$ over elements
$\{n-(j-1)m + 1,\dots,n\}$. Conditional on this, we are interested in the 
action of $\sigma_k$ on $t_j$, i.e. $\{n-jm +1,\dots, n-jm+m\}$. Specifically,
we want to (upper) bound the probability that these elements are mapped to 
themselves. Given $\cap_{\ell=1}^{j-1} \eveF_\ell^c$, each $\sigma_k$ 
will map $\{n-jm +1,\dots, n-jm +m\}$ to one of the ${n-(j-1)m \choose m}$ possibilities with equal probability. 
Further, $\{n-jm +1,\dots, n-jm +m\}$ is not 
a possibility. Therefore, for the purpose of upper bound, we obtain that
\begin{eqnarray}
\Pr\left(\eveF_{j}^c \,\Big|\, \cap_{\ell=1}^{j-1} \eveF_\ell^c\right) 
& \leq &  1 - \left(1 - \frac{1}{{n-(j-1)m \choose m}}\right)^{\sparsity-1} \nonumber \\
& \leq &  1 - \left(1 - \frac{m!}{(n-Lm)^m}\right)^{\sparsity}. \label{eq:c24}
\end{eqnarray}
From \eqref{eq:c21}-\eqref{eq:c24}, we obtain that 
\begin{eqnarray} 
\Pr(\event_1^c) & \leq & \left[1-\left(1-\frac{m!}{(n-Lm)^m}\right)^{\sparsity}\right]^L.\label{eq:c25}
\end{eqnarray}                
Now $Lm = o(n)$ and hence $n-Lm = n(1-o(1))$. Using $1-x = e^{-x} (1+O(x^2))$ 
for small $x \in [0,1)$, we have
\begin{align}
&\left(1-\frac{m!}{(n-Lm)^m}\right) \nonumber\\
=~ & \exp\left(-\frac{m!}{(n-Lm)^m}\right) \left(1 + O\left(\frac{1}{n^{2m}}\right)\right). \label{eq:c26}
\end{align}
Finally, observe that since $\sparsity = O(n^m \log n)$, 
$$ \left(1 + O\left(\frac{1}{n^{2m}}\right)\right)^{\sparsity} = \Theta(1).$$
Thus, from \eqref{eq:c25} and \eqref{eq:c26}, it follows that
\begin{eqnarray} 
\Pr(\event_1^c) & \leq & \left[1- \Theta\left(\exp\left(-\frac{\sparsity m!}{n^m (1-Lm/n)^m}\right)\right) \right]^L \nonumber \\
& \leq & 
\left[1- \Theta\left(\exp\left(-\frac{(1-2\delta)\log n}{(1-n^{-\delta}m)^m}\right)\right) \right]^L \nonumber \\
& \leq & 
\left[1- \Theta\left(\exp\left(-(1-3\delta/2)\log n \right)\right) \right]^L \nonumber \\
& \leq & 
\left[1- \Theta\left(\frac{1}{n^{1-3\delta/2}}\right) \right]^L \nonumber \\
& \leq & \exp\left(-\Omega(Ln^{-1+3\delta/2})\right) \nonumber \\
& \leq & \exp\left(-\Omega(n^{\delta/2})\right) \nonumber \\
& = & O\left(\frac{1}{n^{2m}}\right). \label{eq:c27}
\end{eqnarray}
In above, we have used the fact that $1-x\leq e^{-x}$ for $x \in [0,1]$ and
choice of $L = n^{1-\delta}$. This completes the proof of Theorem \ref{thm:case2}.

\section{Proof of Theorem \ref{thm:case3}: $\lambda_1 = n - n^{\frac{2}{9}-\delta}, \delta > 0$}\label{sec:case3}
So far we have obtained the sharp result that algorithm the sparsest-fit algorithm recovers
$\dist$ up to sparsity essentially $\frac{1}{m!} n^m \log n$ 
for $\lambda$ with $\lambda_1 = n - m$ where $m = O(1)$. Now we
investigate this further when $m$ scales with $n$, i.e. 
$m = \omega(1)$. Let $\lambda_1 = n - \mu$ with $\mu \leq n^{\frac{2}{9}-\delta}$
for some $\delta > 0$. For such $\lambda = (\lambda_1,\dots,\lambda_r)$, 
\begin{eqnarray} 
D_\lambda & = & \frac{n!}{\prod_{i=1}^r \lambda_i !} \nonumber \\
          & \leq & \frac{n!}{\lambda_1 !} \nonumber \\
          & \leq & n^{n-\lambda_1} ~=~n^\mu. \label{eq:30}
\end{eqnarray}
Our interest is in the case when 
$\sparsity \leq (1-\beps) D_\lambda \log \log D_\lambda$
for any $\beps > 0$. For this, the structure of 
arguments will be similar to those used in Theorems 
\ref{thm:case1} and \ref{thm:case2}. Specifically, 
it will be sufficient to establish that 
$\Pr(\event_1^c) = O(1/D_\lambda^2)$, where 
$\event_1$ is the event that permutation 
$\sigma_1 = \id$ satisfies the unique witness property. 

To this end, we order the rows (and corresponding columns) of
the $D_\lambda \times D_\lambda$ matrix $\disthat$ in a
specific manner. Specifically, we are interested in the
$L = 3 n^{\frac{4}{9}-2\delta} \log^3 n$ rows that we call
$t_\ell, 1\leq \ell \leq L$ and they are as follows:
the first row, $t_1$ corresponds to a partition where 
elements $\{1,\dots,\lambda_1\}$ belong to the first 
partition and $\{\lambda_1+1,\dots, n\}$ are partitioned
into remaining $r-1$ parts of size $\lambda_2,\dots,\lambda_r$
in that order.  The partition $t_2$ corresponds to the
one in which the first part contains the $\lambda_1$ elements
$\{1,\dots, n-2\mu,n-\mu+1,\dots,n\}$, while the other $r-1$
parts contain $\{n-2\mu+1,\dots,n-\mu\}$ in that order. 
More generally, for $3 \leq \ell \leq L$, $t_\ell$ contains 
$\{1,\dots,n-\ell \mu,n-(\ell-1)\mu+1,\dots,n\}$ in the first
partition and remaining elements $\{n-\ell \mu+1,\dots, n-(\ell-1)\mu\}$
in the rest of the $r-1$ parts in that order. By our
choice of $L$, $L\mu = o(n)$ and, hence,
the above is well defined.  Next, we bound 
$\Pr(\event_1^c)$ using these $L$ rows. 

Now $\sigma_1 = \id$ and hence $\sigma_1(t_i) = t_i$ 
for all $1\leq i\leq D_\lambda$. Define 
$$\eveF_j = \{\sigma_k(t_j)\neq t_j, \quad \text{for}\quad 2\leq k\leq \sparsity\}, ~~\text{for}~~1\leq j\leq D_\lambda.$$
Then it follows that 
$$ \Pr(\event_1) = \Pr\left(\cup_{j=1}^{D_\lambda} \eveF_j\right).$$
Therefore,  
\begin{eqnarray} 
\Pr(\event_1^c) & = & \Pr\left(\cap_{j=1}^{D_\lambda} \eveF_j^c\right) \nonumber \\
               & \leq & \Pr\left(\cap_{j=1}^L \eveF_j^c\right) \nonumber \\
              & = & \Pr\left(\eveF_1^c\right) \left[\prod_{j=2}^L \Pr\left(\eveF_{j}^c \,\Big|\, \cap_{\ell=1}^{j-1} \eveF_\ell^c\right) \right]. \label{eq:c31}
\end{eqnarray}
First, we bound $\Pr(\eveF_1^c)$. Each permutation $\sigma_k, 1\leq k\leq \sparsity$
maps $t_1$ to one of the $D_\lambda$ possible other $\lambda$ 
partitions with equal probability. Therefore, it follows that  
\begin{eqnarray}
 \Pr\left(\sigma_k(t_1) = t_1\right) & = & \frac{1}{D_\lambda}.\label{eq:c32}
\end{eqnarray}          
Thus,
\begin{eqnarray}
\Pr\left(\eveF_1^c\right) & = & 1 - \Pr\left(\eveF_1\right) \nonumber \\
& = & 1 - \Pr\left(\sigma_k(t_1) \neq t_1, ~2\leq k\leq \sparsity\right) \nonumber \\
& = & 1 - \prod_{k=2}^{\sparsity} \left(1 - \Pr\left(\sigma_k(t_1) = t_1\right)\right) \nonumber \\
& = & 1 - \left(1 - \frac{1}{D_\lambda}\right)^{\sparsity}. \label{eq:c33}
\end{eqnarray}
Next we evaluate $\Pr\left(\eveF_{j}^c \,\Big|\, \cap_{\ell=1}^{j-1} \eveF_\ell^c\right)$
for $2\leq j\leq L$. Given $\cap_{\ell=1}^{j-1} \eveF_\ell^c$, we have (at least partial)
information about the action of $\sigma_k, 2\leq k\leq \sparsity$ over elements
$\{n-(j-1)\mu+1,\dots,n\}$. Conditional on this, we are interested in 
the action of $\sigma_k$ on $t_j$. Given the partial information, each 
of the $\sigma_k$ will map $t_j$ to one of at least $D_{\lambda(j)}$
different options with equal probability for
$\lambda(j) = (\lambda_1 - (j-1)\mu, \lambda_2,\dots,\lambda_r)$ --
this is because the elements $1,\dots, \lambda_1 - (j-1)\mu$ in
the first part and all elements in the remaining $r-1$ parts
are mapped completely randomly conditional on 
$\cap_{\ell=1}^{j-1} \eveF_\ell^c$. Therefore, it follows
that 
\begin{eqnarray}
\Pr\left(\eveF_{j}^c \,\Big|\, \cap_{\ell=1}^{j-1} \eveF_\ell^c\right) 
& \leq &  1 - \left(1 - \frac{1}{D_{\lambda(j)}}\right)^{\sparsity}. \label{eq:c34}
\end{eqnarray}
From \eqref{eq:c31}-\eqref{eq:c34}, we obtain that 
\begin{eqnarray} 
\Pr(\event_1^c) & \leq & \prod_{j=1}^L \left[1-\left(1-\frac{1}{D_{\lambda(j)}}\right)^{\sparsity}\right] \nonumber \\
& \leq & \left[1-\left(1-\frac{1}{D_{\lambda(L)}}\right)^{\sparsity}\right]^L.\label{eq:c35}
\end{eqnarray}           
In above we have used the fact that 
$$ D_{\lambda} ~=~ D_{\lambda(1)} ~\geq~ \dots ~\geq~ D_{\lambda(L)}.$$
Consider
\begin{eqnarray}
\frac{D_{\lambda(j)}}{D_{\lambda(j+1)}} & = & \frac{(n-(j-1)\mu) ! ~(\lambda_1 - j\mu) !}{(n-j\mu) ! ~(\lambda_1 - (j-1)\mu) !} \nonumber \\
& = & \prod_{\ell=0}^{\mu-1}\frac{(n-(j-1)\mu-\ell)}{(\lambda_1 -(j-1)\mu-\ell)} \nonumber \\
& = & \left(\frac{n}{\lambda_1}\right)^\mu \prod_{\ell=0}^{\mu-1}\frac{1-\frac{(j-1)\mu-\ell}{n}}{1-\frac{(j-1)\mu-\ell}{\lambda_1}}
\end{eqnarray}
Therefore, it follows that
\begin{eqnarray}
\frac{D_{\lambda(1)}}{D_{\lambda(L)}} & = & \left(\frac{n}{\lambda_1}\right)^{(L-1)\mu} \prod_{\ell=0}^{(L-1)\mu}\frac{1-\frac{\ell}{n}}{1-\frac{\ell}{\lambda_1}}. \label{eq:c36}
\end{eqnarray}
Using $1+ x \leq e^{x}$ for any $x \in (-1,1)$, $1-x\geq e^{-2x}$ for $x \in (0,1/2)$ and $L\mu = o(n)$, 
we have that for any $\ell, 0\leq \ell \leq (L-1)\mu$
\begin{eqnarray}
\frac{1-\frac{\ell}{n}}{1-\frac{\ell}{\lambda_1}} & = & \frac{1-\frac{\ell}{n} + \frac{\ell}{\lambda_1} - \frac{\ell^2}{n\lambda_1}}{1-\frac{\ell^2}{\lambda^2_1}} \nonumber \\
& \leq & \exp\left(-\frac{\ell^2 - \ell\mu}{n\lambda_1} + \frac{2\ell^2}{\lambda_1^2}\right) \nonumber \\
& \leq & \exp\left(\frac{\ell\mu}{n\lambda_1} + \frac{2\ell^2}{\lambda_1^2}\right). \label{eq:c37}
\end{eqnarray}
Therefore, we obtain
\begin{eqnarray}
\frac{D_{\lambda(1)}}{D_{\lambda(L)}} & \leq & \left(\frac{n}{\lambda_1}\right)^{L\mu} 
\exp\left(\Theta\left(\frac{L^2\mu^3}{n\lambda_1} + \frac{2L^3\mu^3}{\lambda_1^2}\right)\right). \label{eq:38}
\end{eqnarray}
Now 
\begin{eqnarray}
\left(\frac{n}{\lambda_1}\right)^{L\mu} & = & \left(1+\frac{\mu}{\lambda_1}\right)^{L\mu} \nonumber \\
& \leq & \exp\left(\frac{L\mu^2}{\lambda_1}\right). \label{eq:39}
\end{eqnarray}
It can be checked that for given choice of $L, \mu$, we have 
$L\mu^2 = o(\lambda_1)$, $L^3\mu^3 = o(\lambda_1^2)$ and 
$L^2\mu^3 = o(n\lambda_1)$. Therefore, in summary we have that 
\begin{eqnarray}
\frac{D_{\lambda(1)}}{D_{\lambda(L)}} & = & 1 + o(1). \label{eq:310}
\end{eqnarray}
Using similar approximations to evaluate the bound on RHS of \eqref{eq:c35}
along with \eqref{eq:30} yields,
\begin{eqnarray} 
\Pr(\event_1^c) & \leq & \exp\left(-L \exp\left(-\frac{\sparsity}{D_{\lambda(L)}}\right)\right) \nonumber \\
                & = & \exp\left(-L \exp\left(-(1-\beps) \log \log D_\lambda (1+o(1))\right)\right) \nonumber \\
                & \leq & \exp\left(-L\exp\left(-\log \log D_\lambda\right)\right) \nonumber \\
                & = & \exp\left(-\frac{L}{\log D_\lambda}\right) \nonumber \\
                & = & \exp\left(-\frac{3 n^{\frac{4}{9}-2\delta} \log^3 n}{\log D_\lambda}\right) \nonumber \\
                & \leq & \exp\left(-2\log D_\lambda\right) \nonumber \\
                & = & \frac{1}{D_\lambda^2}. \label{eq:312}
\end{eqnarray}
This completes the proof of Theorem \ref{thm:case3}.

\section{Proof of Theorem \ref{thm:case4}: General $\lambda$}\label{sec:case4}
We shall establish the bound on sparsity up to which recovery of $\dist$
is possible from $\disthat$ using the sparsest-fit algorithm for general $\lambda$. 
Let $\lambda = (\lambda_1,\dots,\lambda_r), r\geq 2$ with 
$\lambda_1\geq\dots\geq \lambda_r \geq 1$. As before, let 
\begin{align*}
\sparsity = \|\dist\|_0, ~~ &\Supp{\dist} = \{\sigma_k \in \permgrp :
1\leq k \leq \sparsity\}, \\
\text{and} ~~&\dist(\sigma_k) = p_k, 1\leq k\leq \sparsity.  
\end{align*}
Here $\sigma_k$ and $p_k$ are randomly chosen as per the random
model $\rand$ described in Section \ref{sec:model}. And, we are
given partial information $\disthat$ which is $D_\lambda \times D_\lambda$
matrix with 
$$ D_\lambda = \frac{n!}{\prod_{i=1}^r \lambda_i !}.$$
Finally, recall definition $\alpha = (\alpha_i)_{1\leq i\leq r}$ with
$\alpha_i = \lambda_i/n, 1\leq i\leq r$, 
$$H(\alpha) = -\sum_{i=1}^r \alpha_i \log \alpha_i, \quad\text{and}\quad 
H'(\alpha) = -\sum_{i=2}^r \alpha_i \log \alpha_i.$$
As usual, to establish that the sparsest-fit algorithm recovers $\dist$ from $\disthat$,
we will need to establish ``unique witness'' property as 
``linear independence'' is satisfied due to choice of $p_k$s
as per random model $\rand$. 

For the ease of exposition, we will need an additional notation of
$\lambda$-bipartite graph: it is a complete bipartite graph 
$G^\lambda = (V_1^\lambda \times V_2^\lambda, E^\lambda)$ with
vertices $V_1^\lambda, V_2^\lambda$ having a node each for a
distinct $\lambda$ partition of $n$ and thus 
$|V_1^\lambda| = |V_2^\lambda| = D_\lambda$. Action of a
permutation $\sigma \in \permgrp$, represented by a $0/1$ valued
$D_\lambda \times D_\lambda$ matrix, is equivalent to a perfect
matching in $G^\lambda$. In this notation, a permutation $\sigma$
has ``unique witness'' with respect to a collection of 
permutations, if and only if there is an edge in the matching
corresponding to $\sigma$ that is not present in any other
permutation's matching. 


  Let $\event_L$ denote the event that $L \geq 2$ permutations chosen
  uniformly at random satisfy the ``unique witness'' property. 
  To establish Theorem \ref{thm:case4}, we wish to show that 
  $\Pr(\event^c_\sparsity) = o(1)$ as long as $\sparsity \leq K_1^*(\lambda)$
  where $ K_1^*(\lambda)$ is defined as per \eqref{eq:thmcase4}. 
  To do so, we shall study $\Pr(\event_{L+1}^c \lvert \event_L)$ for $L \geq 1$. 
  Now consider the bipartite graph, $G_L^\lambda$, which is subgraph 
  of $G^\lambda$, formed by the superimposition of the perfect matchings
  corresponding to the $L$ random permutations, $\sigma_i, 1\leq i\leq L$. 
  Now, the probability of $\event_{L+1}^c$ given that $\event_L$ has happened 
  is equal to the probability that a new permutation, generated uniformly
  at random, has its perfect matching so that all its edges end up overlapping 
  with those of $G^\lambda_L$. Therefore, in order to evaluate this probability 
  we count the number such permutations.
  
  For the ease of exposition, we will first count the number of
  such permutations for the cases when $\lambda = (n-1, 1)$ followed
  by $\lambda = (n-2, 2)$. Later, we shall extend the analysis to a 
  general $\lambda$. As mentioned before, for $\lambda = (n-1, 1)$, 
  the corresponding $G^\lambda$ is a complete graph with $n$ nodes
  on left and right. With a bit of abuse of notation, the left and 
  right vertices be labeled $1, 2, \ldots, n$. Now each permutation,
  say $\sigma \in \permgrp$, corresponds to a perfect matching in 
  $G^\lambda$ with an edge from left $i$ to right $j$ if and only if 
  $\sigma(i) = j$. Now, consider $G^\lambda_L$, the superimposition 
  of all the perfect matching of the given $L$ permutations. We want
  to count (or obtain an upper bound on) the number of permutations
  that will have corresponding perfect matching so that all of its
  edges overlap with edges of $G^\lambda_L$. Now each permutation
  maps a vertex on left to a vertex on right. In the graph $G^\lambda_L$, 
  each vertex $i$ on the left has degree of at most $L$. Therefore, 
  if we wish a choose a permutation so that all of its perfect matching's 
  edges overlap with those of $G_\lambda^L$, it has at most $L$ choices
  for each vertex on left. There are $n$ vertices in total on left. Therefore,
  total number of choices are bounded above by $L^n$. From this, we
  conclude that for $\lambda = (n-1,1)$, 
  $$ \Pr(\event_{L+1}^c | \event_L) \leq \frac{L^n}{n!}.$$

  In a similar manner, when $\lambda = (n-2, 2)$, the complete 
  bipartite graph $G^\lambda$ has $D_\lambda = {n \choose 2}$ nodes
  on the left and right; each permutation corresponds to a perfect
  matching in this graph. We label each vertex, on left and right, 
  in $G^\lambda$ by unordered pairs $\set{i,j}$, for $1 \leq i < j \leq n$.
  Again, we wish to bound given $\Pr(\event_{L+1}^c|\event_L)$. For
  this, let $G_L^\lambda$, a subgraph of $G^\lambda$, be obtained by
  the union of edges that belong to the perfect matchings of given $L$
  permutations. We would like to count the number possible permutations
  that will have corresponding matching with edges overlapping with 
  those of $G^\lambda_L$. For this, we consider the $\lfloor n/2\rfloor$ 
  pairs $\set{1,2}, \set{3,4}, \ldots, \{2  \lfloor n/2 \rfloor -1, 2 \lfloor n/2 \rfloor\}$.
  Now if $n$ is even then they end up covering all $n$ elements. If not,
  we consider the last, $n$th element, $\{n\}$ as an additional set. 
  
  Now using a similar argument as before, we conclude that there are at 
  most $L^{\lfloor n/2 \rfloor}$ ways of mapping each of these $\lfloor n/2\rfloor$  
  pairs such that all of these edges overlap with the edges of $G_L^\lambda$. 
  Note that this mapping fixes what each of these $\lceil n/2\rceil$ unordered 
  pairs get mapped to. Given this mapping, there are $2!$ ways of fixing 
  the order in each unordered pair. For example, if an unordered pair
  $\set{i,j}$ maps to unordered pair $\set{k,l}$ there there are $2! = 2$
  options: : $i \mapsto k$, $j \mapsto l$ or $i \mapsto l$, $j \mapsto k$. 
  Thus, once we fix the mapping of each of the $\lceil n/2\rceil$ disjoint 
  unordered pairs, there can be at most $(2!)^{\lceil n/2\rceil}$ permutations 
  with the given mapping of unordered pairs. 
    Finally, note that once the mapping of these $\lfloor n/2\rfloor$ pairs is
    decided, if $n$ is even that there is no element that is left to be 
    mapped. For $n$ odd, since mapping of the $n-1$ elements is decided, so
    is that of $\{n\}$. Therefore, in summary in both even $n$ or odd $n$
    case, there are at most 
  $L^{\lfloor n/2\rfloor} (2!)^{\lceil n/2\rceil}$ permutations that
  have all of the edge of corresponding perfect matching in $G^\lambda$ overlapping
  with the edges of $G^\lambda_L$. Therefore, 
  $$ \Pr(\event_{L+1}^c | \event_L) \leq \frac{L^{\lfloor n/2\rfloor} (2!)^{\lfloor n/2\rfloor}}{n!}.$$

Now consider the case of general $\lambda = (\lambda_1,\lambda_2, \ldots, \lambda_r)$. 
Let $M = \lfloor n/(n-\lambda_1)\rfloor$ and $N = n - M (n-\lambda_1)$. 
Clearly, $0\leq N < n-\lambda_1$. Now we partition the set $\set{1, 2,
    \ldots,n}$ into $M+1$ partitions covering all elements: 
    $\{1,\dots, n-\lambda_1\},\dots, \{(n-\lambda_1)(M-1)+1,\dots, (n-\lambda_1)M\}$
    and $\{(n-\lambda_1)M+1,\dots, n\}$. 
    As before, for the purpose of upper bounding the number of permutations
    that have corresponding perfect matchings in $G^\lambda$ overlapping
    with edges of $G_L^\lambda$, each of the first $M$ partitions can be mapped 
    in $L$ different ways; in total at most $L^M$ ways. For each of these mappings, 
    we have options at the most 
    $$(\lambda_2! \lambda_3! \ldots \lambda_r!)^M.$$ 
    Given the mapping of the first $M$ partitions, the mapping of the $N$ 
    elements of the $M+1$st partition is determined (without ordering). Therefore,
    the additional choice is at most $N!$.  
    In summary, the total number of permutations can be at most
$$   L^M \left(\prod_{i=2}^r \lambda_i ! \right)^M N!.$$
Using this bound, we obtain
    \begin{eqnarray}
\Pr\left(\event_{L+1}^c \lvert \event_L\right) & \leq  & \frac{1}{n!}
L^M \left(\prod_{i=2}^r \lambda_i ! \right)^M N!. \label{eq:c41}
\end{eqnarray}
Let, 
$$x_L \stackrel{\triangle}{=}\frac{1}{n!}
L^M \left(\prod_{i=2}^r \lambda_i ! \right)^M N!.$$
Note that $\event_{k+1} \subset \event_k$ for $k\geq 1$. Therefore,
it follows that 
\begin{eqnarray}
\Pr\left(\event_\sparsity\right) & = & \Pr\left(\event_\sparsity \cap \event_{\sparsity-1}\right) \nonumber \\
& = & \Pr\left(\event_\sparsity \lvert \event_{\sparsity-1}\right) \Pr\left(\event_{\sparsity-1}\right).\label{eq:42}
\end{eqnarray}
Recursive application of argument behind \eqref{eq:42} and fact that $\Pr(\event_1) = 1$, we
have
\begin{eqnarray}
\Pr\left(\event_\sparsity\right) & = & \Pr\left(\event_1\right)\prod_{L=1}^{\sparsity-1} \Pr\left(\event_{L+1} \lvert \event_L\right) \nonumber \\
& = & \prod_{L=1}^{\sparsity-1} \left(1-\Pr\left(\event_{L+1}^c \lvert \event_L\right)\right) \nonumber \\
& = & \prod_{L=1}^{\sparsity-1} \left(1-x_{L}\right) \nonumber \\
& \geq & 1 - \left(\sum_{L=1}^{\sparsity-1} x_L \right). \label{eq:c43} 
\end{eqnarray}
Using \eqref{eq:c41}, it follows that $x_{k+1} \geq x_k$ for $k \geq 1$. 
Therefore,
\begin{eqnarray}
\sum_{L=2}^\sparsity x_L & \leq & \sparsity x_\sparsity \nonumber \\
& \leq & \frac{1}{n!}\sparsity^{M+1} \left(\prod_{i=2}^r \lambda_i ! \right)^M N! \nonumber\\
&=& \frac{1}{n!} \sparsity^{M+1} \left( \frac{n!}{\lambda_1! D_\lambda} \right)^M N! \nonumber \\
&=& \frac{\sparsity^{M+1}}{D_\lambda^M} \left( \frac{n!}{\lambda_1!}\right)^M \frac{N!}{n!} \nonumber \\
&=& \frac{\sparsity^{M+1}}{D_\lambda^M}  \left(\frac{n!}{\lambda_1! (n-\lambda_1)!}\right)^M \frac{N! ((n-\lambda_1)!)^M}{n!}.\label{eq:c45}
\end{eqnarray}
Since $n = N + M(n - \lambda_1)$, we have a binomial and a multinomial coefficient in RHS of \eqref{eq:c45}. We simplify this expression by obtaining an approximation for a multinomial coefficient through Stirling's approximation. For that, first consider a general multinomial coefficient $m!/(k_1! k_2! \ldots k_l!)$ with $m = \sum_i k_i$. Then, using the Stirling's approximation $\log n! = n \log n -n + 0.5 \log n + O(1)$, for any $n$, we obtain 

\begin{align*}
&\log \left( \frac{m!}{k_1! k_2! \ldots k_l!} \right) \nonumber\\
 =~ &m \log m - m + 0.5 \log m + O(1) - \\
&\sum_{i = 1}^l \left(k_i \log k_i - k_i + 0.5 \log k_i + O(1)\right) \nonumber \\
 = ~&m \sum_{i = 1}^l \frac{k_i}{m} \log \frac{m}{k_i} + 0.5 \log \frac{m}{k_1 k_2 \ldots k_l} - O(l) 
\end{align*}

Thus, we can write
\begin{align}
  &M \log \frac{n!}{\lambda_1! (n-\lambda_1)!} \nonumber\\
\;\; = \;\; &Mn  \alpha_1 \log\frac{1}{\alpha_1} +  M n (1 - \alpha_1) \log \frac{1}{1 - \alpha_1} \label{eq:j11}\\
 &+ 0.5 \log \frac{1}{n^M \alpha_1^M (1 - \alpha_1)^M} - O(M) \nonumber
\end{align}
where $\alpha_1 = \lambda_1/n$. Similarly, we can write 
\begin{align}
  &\log \frac{n!}{N! ((n - \lambda_1)!)^M} \nonumber\\
\;\;= \;\;&n \delta \log \frac{1}{\delta} + M n (1 - \alpha_1) \log \frac{1}{1 - \alpha_1} \label{eq:j12}\\
&+ 0.5 \log \frac{1}{n^M \delta ( 1 - \alpha_1)^M} - O(M) \nonumber  
\end{align}
where $\delta = N/n$. It now follows from \eqref{eq:j11} and \eqref{eq:j12} that 
\begin{align}
   &M \log \frac{n!}{\lambda_1! (n - \lambda_1)!} - \log \frac{n!}{N! ((n - \lambda_1)!)^M} \nonumber \\
= ~&- M n \alpha_1 \log \alpha_1 + \delta n \log \delta \label{eq:s3} \\
&+ 0.5 \log \frac{\delta}{\alpha_1^M} + O(M) \nonumber
\end{align}
Since $\delta < 1$, $\delta n \log \delta \leq 0$ and $\log (\delta/ \alpha_1^M) \leq - M \log \alpha_1$. Thus, we can write 
\begin{align}
  &M \log \frac{n!}{\lambda_1! (n - \lambda_1)!} - \log \frac{n!}{N! ((n - \lambda_1)!)^M} \nonumber \\
\leq ~& M n \alpha_1 \log (1/\alpha_1) + O( M \log (1/\alpha_1)) \nonumber \\
=~ &O(M n \alpha_1 \log (1/\alpha_1)) \label{eq:s4}
\end{align}

It now follows from \eqref{eq:c45}, \eqref{eq:s3} and \eqref{eq:s4} that 
\begin{align}
&\log \left( \sum_{L = 2}^K x_L \right) \nonumber \\
\leq &(M+1) \log K - M \log D_\lambda + O(M n \alpha_1 \log (1/\alpha_1)) \label{eq:s5}
\end{align}
Therefore, for $\Pr(\event_\sparsity) = 1 - o(1)$, a sufficient condition is 
\begin{align}
&\log K + \frac{c\log n}{M+1} \nonumber \\
\leq & \frac{M}{M+1} \log D_\lambda - \frac{M}{M+1}O(n \alpha_1 \log (1/\alpha_1)) \label{eq:s6}
\end{align}
for some $c > 0$. We now claim that $\log n = O(M n \alpha_1 \log(1/\alpha_1))$. The claim is clearly true for $\alpha_1 \to \theta$ for some $0 < \theta < 1$. Now suppose $\alpha_1 \to 1$. Then, $M \geq 1/(1-\alpha_1) - 1 = \alpha_1/(1 - \alpha_1) = x$, say. This implies that $M \alpha_1 \log(1/\alpha_1) \geq  \alpha_1 x \log( 1 + 1/x) \to 1$ as $\alpha_1 \to 1$. Thus, $M n \alpha_1 \log(1/\alpha_1) = n(1 + o(1))$ for $\alpha_1 \to 1$ as $n \to \infty$. Hence, the claim is true for $\alpha_1 \to 1$ as $n \to \infty$. Finally, consider $\alpha_1 \to 0$ as $n \to \infty$. Note that the function $h(x) = x \log(1/x)$ is increasing on $(0, \epsilon)$ for some $0 < \epsilon < 1$. Thus, for $n$ large enough, $n \alpha_1 \log(1/\alpha_1) \geq \log n$ since $\alpha_1 \geq 1/n$. Since $M \geq 1$, it now follows that $M n \alpha_1 \log (1/\alpha_1) \geq \log n$ for $n$ large enough and $\alpha_1 \to 0$. This establishes the claim. 



Since $\log n = O(M n \alpha_1 \log(1/\alpha_1))$, it now follows that \eqref{eq:s6} is implied by 
\begin{align}
\log K & \leq  \frac{M}{M+1} \log D_\lambda - \frac{M}{M+1}O(n \alpha_1 \log (1/\alpha_1)) \nonumber\\ 
& =  \frac{M}{M+1} \log D_\lambda \left[ 1 - \frac{O(n\alpha_1 \log (1/\alpha_1))}{\log D_\lambda} \right]\label{eq:sjj1}
\end{align}

Now consider $D_\lambda = n!/ (\lambda_1! \lambda_2! \ldots \lambda_r!)$. Then, we claim that for large $n$
\begin{equation}
\log D_\lambda \geq  0.5 n H(\alpha). \label{eq:c414}  
\end{equation}
In order to see why the claim is true, note that Stirling's approximation suggests, 
\begin{eqnarray*}
\log n! & = & n \log n - n + 0.5 \log n + O(1) , \nonumber \\
\log \lambda_i ! & = & \lambda_i\log \lambda_i - \lambda_i + 0.5 \log \lambda_i + O(1). \label{eq:c414x}
\end{eqnarray*}
Therefore, 
\begin{eqnarray*}
\log D_\lambda & \geq & n H(\alpha) + 0.5 \log(n/\lambda_1) - \sum_{i = 2}^r 0.5 (O(1) + \log \lambda_i). 
\end{eqnarray*}
Now consider, 
\begin{align} 
  &\lambda_i \log(n/\lambda_i) - \log \lambda_i - O(1) \nonumber \\
= &\left( \lambda_i - \frac{\log \lambda_i}{\log(n/\lambda_i)} \right)
\log(n/\lambda_i) - O(1) \label{eq:js11}
\end{align}
Since $\lambda_i \leq n/2$ for $i \geq 2$, $\log(n/\lambda_i) \geq \log 2$. Thus, the first term in the RHS of \eqref{eq:js11} is non-negative for any $\lambda_i \geq 1$. In addition, for every $\lambda_i$, either $\lambda_i - \log \lambda_i \to \infty$ or $\log(n/\lambda_i) \to \infty$ as $n \to \infty$. Therefore, the term on the RHS of \eqref{eq:js11} is asymptotically non-negative. Hence,
\begin{eqnarray}
\log D_\lambda & \geq & 0.5 n H(\alpha). 
\end{eqnarray}
Thus, it now follows from \eqref{eq:c414} that \eqref{eq:sjj1} is implied by 
\begin{eqnarray*}
\qquad \log K & \leq & \frac{M}{M+1} \log D_\lambda \left[ 1 - \frac{O(\alpha_1 \log (1/\alpha_1))}{H(\alpha)} \right].\label{eq:c415}
\end{eqnarray*} 
That is, we have ``unique witness'' property satisfied as long as 
\begin{eqnarray}
\sparsity & = & O\left(D_\lambda^{\gamma(\alpha)}\right), 
\end{eqnarray}
where
\begin{eqnarray} 
\gamma(\alpha) & = & \frac{M}{M+1} \left[ 1 - C' \frac{H(\alpha) -
    H'(\alpha)} {H(\alpha)} \right], 
\end{eqnarray}
and $C'$ is some constant.
This completes the proof of Theorem \ref{thm:case4}. 

\section{Proof of Theorem \ref{thm:converse}: Limitation on Recovery}\label{sec:thmconverse} 

In order to make a statement about the inability of {\em any} algorithm to recover $\dist$ using $\disthat$, we rely on the formalism of classical information theory. In particular, we establish a bound on the sparsity of $\dist$ beyond which recovery is not asymptotically reliable (precise definition of asymptotic reliability is provided later). 


\subsection{Information theory preliminaries}\label{sub:ITprelims} Here we recall some
necessary Information Theory preliminaries. Further details can
be found in the book by Cover and Thomas~\cite{CT06}. 

Consider a discrete random variable $X$ that is uniformly distributed over
a finite set $\crX$. Let $X$ be {\em transmitted} over a 
{\em noisy} channel to a receiver; suppose the receiver receives a
random variable $Y$, which takes values in a finite set $\crY$.  
Essentially, such ``transmission over noisy channel'' setup describes 
any two random variables $X, Y$ defined through a joint
probability distribution over a common probability space. 

Now let $\hat{X} = g(Y)$ be an estimation of the transmitted
information that the receiver produces based on the observation $Y$ using some
function $g : \crY \to \crX$. Define probability of error 
as $p_{\text{err}} = \Pr(X \neq \hat{X})$. Since $X$ is
uniformly distributed over $\crX$, it follows that
\begin{equation}
  \label{eq:104}
  p_{\text{err}} = \frac{1}{\abs{\crX}} \sum_{x \in \crX} \Pr(g(Y)
  \neq x \lvert x).
\end{equation}
Recovery of $X$ is called asymptotically reliable if $p_{\text{err}}
\to 0$ as $\abs{\crX} \to \infty$. Therefore, in order to show that recovery is not asymptotically reliable, it is sufficient to prove that $p_{\text{err}}$ is bounded away from $0$ as $\abs{\crX} \to \infty$. In order to obtain a lower bound on $p_{\text{err}}$, we use Fano's inequality:
\begin{eqnarray}
H(X|\hat{X}) & \leq & 1 + p_{\text{err}} \log \abs{\crX}. \label{eq:c101}
\end{eqnarray}
Using \eqref{eq:c101}, we can write
\begin{eqnarray}
H(X) & = & I(X; \hat{X}) + H(X|\hat{X}) \nonumber \\
     & \leq & I(X; \hat{X}) + p_{\text{err}} \log \abs{\crX} + 1 \nonumber \\
     & \stackrel{(a)}{\leq} & I(X; Y) + p_{\text{err}} \log \abs{\crX} + 1 \nonumber \\ 
     & = & H(Y) - H(Y | X) + p_{\text{err}} \log \abs{\crX} + 1 \nonumber \\
     & \leq & H(Y) + p_{\text{err}} \log \abs{\crX} + 1, \label{eq:c102}
\end{eqnarray}
where we used $H(Y|X) \geq 0$ for a discrete\footnote{The counterpart of
this inequality for a continuous valued random variable is not true. 
This led us to study the limitation of recovery algorithm 
over model $\rT$ rather than $\rand$.} 
valued random variable. The inequality 
(a) follows from the data processing inequality: if we have Markov
chain $X \to  Y \to \hat{X}$, then $I(X; \hat{X}) \leq I(X; Y)$. Since $H(X) = \log \abs{\crX}$, 
from \eqref{eq:c102} we obtain
\begin{eqnarray}
p_{\text{err}} & \geq & 1- \frac{H(Y) + 1}{\log \abs{\crX}}. \label{eq:c103}
\end{eqnarray}
Therefore, to establish that probability of error is bounded away from zero, it is
sufficient to show that 
\begin{eqnarray}
\frac{H(Y) + 1}{\log \abs{\crX}} & \leq & 1-\delta, \label{eq:c104}
\end{eqnarray}
for any fixed constant $\delta > 0$. 

\subsection{Proof of theorem \ref{thm:converse}.} 

Our goal is to show that when $\sparsity$ is large enough (in particular, as claimed in the statement of Theorem \ref{thm:converse}), the probability of error of {\em any} recovery algorithm is uniformly bounded away from $0$. For that, we first fix a recovery algorithm, and then utilize the above setup to show that recovery is not asymptotically reliable when $\sparsity$ is large. Specifically, we use \eqref{eq:c104}, for which we need to identify random variables $X$ and $Y$.

To this end, for a given $\sparsity$ and $T$, let $\dist$ be generated as per the random model $\rT$. Let
random variable $X$ represent the support of
function $\dist$ i.e., $X$ takes values 
in $\crX = \permgrp^\sparsity$. Given $\lambda$,
let $\disthat$ be the partial information that the recovery 
algorithm uses to recover $\dist$. Let random
variable $Y$ represent $\disthat$, the $D_\lambda\times D_\lambda$
matrix. Let $h = h(Y)$ denote the estimate of $\dist$, and 
$g=g(Y) = \supp{h}$ denote the estimate of the 
support of $\dist$ produced by the given recovery algorithm. Then,
\begin{eqnarray}
\Pr\left(h \neq \dist\right) & \geq & \Pr\left(\supp{h} \neq \supp{\dist}\right) \nonumber \\
& = & \Pr\left(g(Y) \neq X\right). \label{eq:c51}
\end{eqnarray}
Therefore, in order to uniformly lower bound the probability of error of the recovery algorithm, it is sufficient to lower bound its probability of making an error in recovering the support of $\dist$. Therefore, we focus on
$$ p_{\text{err}} = \Pr\left(g(Y) \neq X\right).$$

It follows from the discussion in Section \ref{sub:ITprelims} that in order to show that $p_{\text{err}}$ is uniformly bounded away from $0$, it is sufficient to show that for some constant $\delta > 0$ \begin{eqnarray} \frac{H(Y)+1}{\log \abs{\crX}} & \leq & 1 - \delta. \label{eq:c52} \end{eqnarray} Observe that $\abs{\crX} = (n!)^{\sparsity}$. Therefore, using Stirling's approximation, it follows that \begin{eqnarray} \log \abs{\crX} & = & (1+o(1)) \sparsity n \log n. \label{eq:c53} \end{eqnarray} Now $Y = \disthat$ is a $D_\lambda\times D_\lambda$ matrix.  Let $Y = [Y_{ij}]$ with $Y_{ij}, 1\leq i, j \leq D_{\lambda}$, taking values in $\{1,\dots, \sparsity T\}$; it is easy to see that $H(Y_{ij}) \leq \log \sparsity T$. Therefore, it follows that 
\begin{eqnarray}
  H(Y) & \leq & \sum_{i, j = 1}^{D_\lambda} H(Y_{ij}) \nonumber \\
     & \leq & D_{\lambda}^2 \log \sparsity T ~=~ D_{\lambda}^2 \left(\log \sparsity + \log T\right). \label{eq:c54}
\end{eqnarray}
For small enough constant $\delta > 0$, it is easy to see that  
the condition of \eqref{eq:c52} will follow if $\sparsity$ satisfies the following two inequalities:
\begin{eqnarray}
\frac{D_{\lambda}^2 \log \sparsity}{\sparsity n \log n} \leq  \frac{1}{3}\left(1+\delta\right) 
\quad \Leftarrow\quad  \frac{K}{\log K} \geq \frac{3 (1-\delta/2)D_{\lambda}^2}{ n \log n}, \label{eq:c55a} \\
\frac{D_{\lambda}^2 \log T}{\sparsity n \log n}  \leq \frac{1}{3}\left(1+\delta\right)
\quad\Leftarrow\quad K \geq \frac{3 (1-\delta/2) D_{\lambda}^2\log T}{ n \log n}. \label{eq:c55b}
\end{eqnarray}
In order to obtain a bound on $\sparsity$ from \eqref{eq:c55a}, consider the following: for large numbers $x, y$, 
let $y = (c+\beps) x \log x$, for some constants $c, \beps > 0$. Then, $\log y = \log x + \log \log x + \log (c+\beps)$
which is $(1+o(1)) \log x$. Therefore, 
\begin{eqnarray}
\frac{y}{\log y} & = & \frac{c+\beps}{1+o(1)} x ~\geq~cx, \label{eq:c55a1}
\end{eqnarray}
for $x \to \infty$ and constants $c,\beps > 0$. Also, observe that $y/ \log y$ is a non-decreasing function; hence, it follows that for $y \geq (c+\beps) x \log x$, $y / \log y \geq cx$ for large $x$. Now take $x = \frac{D_\lambda^2}{n\log n}$, $c = 3$, $\beps = 1$ and $y = K$. Note that $ D_\lambda \geq n$ for all $\lambda$ of interest; therefore, $x \to \infty$ as $n\to \infty$. 
Hence, \eqref{eq:c55a} is satisfied for the choice of 
\begin{eqnarray}
\sparsity & \geq & \frac{4 D_{\lambda}^2}{ n \log n} \left(\log \frac{D_{\lambda}^2}{ n \log n} \right). \label{eq:c55a2}
\end{eqnarray}
From \eqref{eq:c52}, \eqref{eq:c55a}, \eqref{eq:c55b}, and \eqref{eq:c55a2}
 it follows that the probability of error of
any algorithm is at least $\delta > 0$ for $n$ large enough and
any $\lambda$ if 
\begin{eqnarray}
\sparsity & \geq & \frac{4 D_{\lambda}^2}{ n \log n} \left[ \log \left(\frac{D_{\lambda}^2}{ n \log n} \vee T\right) \right]. \end{eqnarray}
This completes the proof of Theorem \ref{thm:converse}.

\section{Conclusion} \label{sec:conclusion}

In summary, we considered the problem of {\em exactly} recovering a non-negative function over the space of permutations from a given partial set of Fourier coefficients. This problem is motivated by the wide ranging applications it has across several disciplines. This problem has been widely studied in the context of discrete-time functions in the recently popular {\em compressive sensing} literature. However, unlike our setup, where we want to perform exact recovery from a {\em given set} of Fourier coefficients, the work in the existing literature pertains to the choice of a limited set of Fourier coefficients that can be used to perform exact recovery. 

Inspired by the work of Donoho and Stark~\cite{DS89} in the context of discrete-time functions, we focused on the recovery of non-negative functions with a sparse support (support size $\ll$ domain size). Our recovery scheme consisted of finding the function with the sparsest support, consistent with the given information, through $\ell_0$ optimization. As we showed through some counter-examples, this procedure, however, will not recover the exact solution in all the cases. Thus, we identified sufficient conditions under which a function can be recovered through $\ell_0$ optimization. For each kind of partial information, we then quantified the sufficient conditions in terms of the ``complexity'' of the functions that can be recovered. Since the sparsity (support size) of a function is a natural measure of its complexity, we quantified the sufficient conditions in terms of the sparsity of the function. In particular, we proposed a natural random generative model for the functions of a given sparsity. Then, we derived bounds on sparsity for which a function generated according to the random model satisfies the sufficient conditions with a high probability as $n \to \infty$. Specifically, we showed that, for partial information corresponding to partition $\lambda$, the sparsity bound essentially scales as $D_\lambda^{M/(M+1)}$. For $\lambda_1/n \to 1$, this bound essentially becomes $D_\lambda$ and for $\lambda_1/n \to 0$, the bound essentially becomes $D_\lambda^{1/2}$.

Even though we found sufficient conditions for the recoverability of functions by finding the sparsest solution, $\ell_0$ optimization is in general computationally hard to carry out. This problem is typically overcome by considering its convex relaxation, the $\ell_1$ optimization problem. However, we showed that $\ell_1$ optimization fails to recover a function generated by the random model with a high probability. Thus, we proposed a novel iterative algorithm to perform $\ell_0$ optimization for functions that satisfy the sufficient conditions, and extended it to the general case when the underlying distribution may not satisfy the sufficient conditions and the observations maybe noisy.

We studied the limitation of any recovery algorithm by means of information theoretic tools. While the bounds we obtained are useful in general, due to technical limitations, they do not apply to the random model we considered. Closing this gap and understanding recovery conditions in the presence of noise are natural next steps.

\bibliographystyle{IEEEtran}
\bibliography{bibdb}

\appendix[Proof of Auxiliary Lemma] \label{ap1}

Here we present the proof of Lemma \ref{lem:analysis}. For this, first consider the limit $\alpha_1 \uparrow 1$. Specifically, let $\alpha_1 = 1-\beps$, for a very small positive $\beps$. Then, $\sum_{i=2}^r \alpha_i = 1-\alpha_1 = \beps$. By definition, we have $H'(\alpha)/H(\alpha) \leq 1$; therefore, in order to prove that $H'(\alpha)/H(\alpha) \to 1$ as $\alpha_1 \uparrow 1$, it is sufficient to prove that $H'(\alpha)/H(\alpha) \geq 1 - o(1)$ as $\alpha_1 \uparrow 1$. For that, consider 
\begin{eqnarray} \frac{H'(\alpha)}{H(\alpha)} & = &
  \frac{H'(\alpha)}{ \alpha_1 \log (1/\alpha_1) + H'(\alpha) } \nonumber \\
  & = & 1 - \frac{\alpha_1\log(1/\alpha_1)}{\alpha_1 \log (1/\alpha_1) + H'(\alpha) }.\label{eq:apx1} 
\end{eqnarray} 
In order to obtain a lower bound, we minimize $H'(\alpha)/H(\alpha)$ over $\alpha \geq 0$. It follows from \eqref{eq:apx1} that, for a given $\alpha_1 = 1-\beps$, $H'(\alpha)/H(\alpha)$ is minimized for the choice of $\alpha_i, i\geq 2$ that minimizes $H'(\alpha)$. Thus, we maximize $\sum_{i=2}^r \alpha_i\log \alpha_i$ subject to $\alpha_i \geq 0$ and $\sum_{i=2}^r \alpha_i = 1-\alpha_1 = \beps$. Here we are maximizing a convex function over a convex set. Therefore, maximization is achieved on the boundary of the convex set. That is, the maximum is $\beps \log \beps$; consequently, the minimum value of $H'(\alpha)= \beps \log (1/\beps)$.  Therefore, it follows that for $\alpha_1 = 1-\beps$, 
\begin{eqnarray}
  1 ~\geq~ \frac{H'(\alpha)}{H(\alpha)} & \geq & 1 - \frac{-(1-\beps) \log (1-\beps) }{ \beps\log (1/\beps) -(1-\beps) \log (1-\beps) } \nonumber \\
  & \approx & 1 - \frac{\beps}{\beps \log (1/\beps) + \beps} \nonumber \\
  & \approx & 1 - \frac{1}{1+ \log (1/\beps)} \nonumber \\
  & \stackrel{\beps\to 0}{\rightarrow} & 1. \label{eq:apx2} 
\end{eqnarray}

To prove a similar claim for $\alpha_1\downarrow 0$, let $\alpha_1 =\beps$
for a small, positive $\beps$. Then, it follows that $r = \Omega(1/\beps)$
since $\sum_{i=1}^r \alpha_i = 1$ and $\alpha_1 \geq \alpha_i$ for all 
$i, 2\leq i\leq r$. Using a convex maximization based
argument similar to the one we used above, it can be checked that $H'(\alpha) = \Omega(\log (1/\beps))$. 
Therefore, it follows that $\alpha_1 \log(1/\alpha_1)/H'(\alpha) \to 0$
as $\alpha_1 \downarrow 0$.  That is, $H'(\alpha)/H(\alpha) \to 1$ as
$\alpha_1 \downarrow 0$. This completes the proof of Lemma~\ref{lem:analysis}.

\begin{IEEEbiographynophoto}{Srikanth Jagabathula} 
Srikanth Jagabathula received the BTech degree in Electrical
Engineering from the Indian Institute of Technology (IIT) Bombay in
2006, and the MS degree in Electrical Engineering and Computer Science
from the Massachusetts Institute of Technology (MIT), Cambridge, MA in
2008. He is currently a doctoral student in the Department of
Electrical Engineering and Computer Science at MIT. His research
interests are in the areas of revenue management, choice modeling,
queuing systems, and compressed sensing. He received the President of
India Gold Medal from IIT Bombay in 2006. He was also awarded the
"Best Student Paper Award" at NIPS 2008 conference, the Ernst
Guillemin award for the best EE SM Thesis, first place in the MSOM
student paper competition in 2010. \end{IEEEbiographynophoto}

\begin{IEEEbiographynophoto}{Devavrat Shah} Devavrat Shah is currently
  a Jamieson career development associate professor with the
  department of electrical engineering and computer science, MIT. He
  is a member of the Laboratory of Information and Decision Systems
  (LIDS) and affiliated with the Operations Research Center (ORC). His
  research focus is on theory of large complex networks which includes
  network algorithms, stochastic networks, network information theory
  and large scale statistical inference.  He received his BTech degree
  in Computer Science and Engg. from IIT-Bombay in 1999 with the honor
  of the President of India Gold Medal. He received his Ph.D. from the
  Computer Science department, Stanford University in October 2004. He
  was a post-doc in the Statistics department at Stanford in 2004-05.
  He was co-awarded the best paper awards at the IEEE INFOCOM '04, ACM
  SIGMETRICS/Performance '06; and best student paper awards at Neural
  Information Processing Systems '08 and ACM SIGMETRICS/Performance
  '09. He received 2005 George B. Dantzig best dissertation award from
  the INFORMS. He received the first ACM SIGMETRICS Rising Star Award
  2008 for his work on network scheduling
  algorithms.  \end{IEEEbiographynophoto}

\ifCLASSOPTIONcaptionsoff
  \newpage
\fi

\end{document}